\documentclass[10pt,a4paper]{article}

\usepackage{amsmath}
\usepackage{graphicx,latexsym,euscript,makeidx,color}
\usepackage{amsmath,amsfonts,amssymb,amsthm,mathrsfs,dsfont}
\usepackage{enumerate}

\usepackage{amssymb}
\usepackage{amsmath}
\usepackage{amsfonts}
\usepackage{graphicx}
\usepackage{graphicx}          
\usepackage{color}

\usepackage[colorlinks,linkcolor=blue,anchorcolor=green,citecolor=red]{hyperref}

\usepackage{geometry}
\font\tenbb=msbm10 \font\sevenbb=msbm7 \font\fivebb=msbm5
\newfam\bbfam
\scriptscriptfont\bbfam=\fivebb \textfont\bbfam=\tenbb
\scriptfont\bbfam=\sevenbb

\newtheorem{theorem}{\indent Theorem}[section]
\newtheorem{definition}[theorem]{\indent Definition}
\newtheorem{proposition}[theorem]{\indent Proposition}

\newtheorem{lemma}[theorem]{\indent Lemma}
\newtheorem{remark}[theorem]{\indent Remark}
\newtheorem{example}[theorem]{\indent Example}

\makeatletter
   
   \@addtoreset{equation}{section}
\makeatother

\allowdisplaybreaks[4]

\begin{document}

\title{\bf Deterministic Dynamic Stackelberg Games: Time-Consistent Open-Loop Solution\thanks{This work is supported in part by the National Natural Science Foundation of China (61773222, 11871369, 61973172, 62173191). 
}
}
\author{Yuan-Hua Ni\thanks{College of Artificial Intelligence, Nankai University, Tianjin 300350, P.R. China. Email: {\tt yhni@nankai.edu.cn}. }~~~~~Liping Liu\thanks{College of Artificial Intelligence, Nankai University, Tianjin 300350, P.R. China. Email: {\tt liulipingnk@163.com}.}~~~~~Xinzhen Zhang\thanks{School of Mathematics, Tianjin University, Tianjin 300352, P.R. China. Email: {\tt xzzhang@tju.edu.cn}.}}
\maketitle

{\bf Abstract:} In this paper, the known deterministic linear-quadratic Stackelberg game is revisited, whose open-loop Stackelberg solution actually possesses the nature of time inconsistency. To handle this time inconsistency, {a two-tier game framework is introduced, where the upper-tier game works according to Stackelberg's scenario with a leader and a follower, and two lower-tier intertemporal games give the follower's and leader's equilibrium response mappings that mimic the notion of time-consistent open-loop equilibrium control in existing literature.
The resulting open-loop equilibrium solution of the two-tier game} is shown to be weakly time-consistent in the sense that the adopted policies will no longer be denied in the future only if past policies are consistent with the equilibrium policies. On the existence and uniqueness of such a solution, necessary and sufficient conditions are obtained, which are characterized via the solutions of several Riccati-like equations.



\section{Introduction}


In its original setting, Stackelberg game is static  \cite{Stackelberg}, which is firstly formulated by H. von Stackelberg in 1934 to meet the markets with some firms having dominating power over others;
within the context of two-person nonzero-sum games, the dominator first announces her policy who is called the leader, then the other player, called the follower, reacts to minimize his cost functional, and finally the leader makes her optimal decision under the follower's best response. Then, the leader's optimal policy and the follower's best response form the known Stackelberg equilibrium or Stackelberg solution.

The extension of Stackelberg games to their dynamic setting is traced back to the early 1970s  \cite{Simaan-1,Simaan-2}, where both the local and global Stackelberg solutions are introduced. Local Stackelberg solution corresponds to the case  where the leader announces her policy to the follower stagewisely, namely, the leader has stagewise dominating power over the follower. Local Stackelberg solution is seeked via a backward recursion as that in dynamic programming, and at each step of the recursion we are facing a static Stackelberg game with the leader's decision information being the state variable at that time.
Clearly, such local Stackelberg solution is of feedback nature, which is also called feedback Stackelberg solution.
Another important solution concept of dynamic Stackelberg game is the global Stackelberg solution, where the leader has dominating power over the follower on the lifetime horizon, i.e., the leader looks at the time horizon as a whole and announces her policy over the horizon before the game starts. Corresponding to the underlining information structure (the whole of decision information sets), global Stackelberg solutions are classified into several types, and open-loop Stackelberg solution is the representative one that is firstly studied in \cite{Simaan-1} and corresponds to the open-loop information structure. For the overview of global Stackelberg solutions, we are referred to the monograph \cite{Basar} and recent work \cite{Bensoussan}. 

Open-loop solution in dynamic games has a long history, which may date back early to the work \cite{Berkovitz}, and open-loop Stackelberg solution has gained much attention during the last half century; to name a few of the literature, see \cite{Freiling,Kydland-2,Shi-1,Juan-1,Yong-3} besides \cite{Basar,Bensoussan,Simaan-1}.
In \cite{Simaan-1}, linear-quadratic (LQ, for short) dynamic Stackelberg games are formulated and solved in a Hilbert space setting, and sufficient condition on the existence and uniqueness of open-loop Stackelberg solution is presented.
In \cite{Kydland-2}, dynamic Stackelberg games with special structure (called dominant-player game there) are investigated, and the distinctions between open-loop solution, feedback solution and closed-loop solution are clarified from the viewpoint of economist.
Sufficient conditions on the existence of open-loop solution is presented in \cite{Freiling} for LQ differential games via Lyapunov-type approach. In \cite{Yong-3}, a leader-follower stochastic LQ differential game is investigated, where the coefficients of the controlled system and cost functionals are random and the weighting matrices are allowed to be indefinite. If two stochastic Riccati equations are solvable, the open-loop Stackelberg solution is shown to admit feedback representation \cite{Yong-3}.
Necessary and sufficient conditions are presented in \cite{Juan-1} to characterize the existence and uniqueness of open-loop solution for a deterministic LQ dynamic Stackelberg game, which resort to two discrete-time Riccati equations. A Stackelberg stochastic differential game with asymmetric information is studied in \cite{Shi-1}, which makes use of the stochastic maximum principle and verification theorem with partial information to derive the open-loop equilibrium solution.

In the survey paper \cite{Sethi}, the authors mention through a few sentences that open-loop Stackelberg solution is not time-consistent in general, namely, ``This means that given an opportunity to revise her strategy at any future time after the initial time, the leader would benefit by choosing another strategy than the one she chose at the initial time."  Except for two literal paragraphs, yet no more information has been provided in \cite{Sethi} about the time inconsistency of open-loop Stackelberg solution.
To the best of the authors, the first that proposes the time inconsistency of open-loop Stackelberg solution is the work \cite{Kydland-2}, which {checks} a special case where the cost functional can be written in terms of the decision variables (control inputs) only. By computing the first-order conditions for an optimum of the leader, one has the difference of these first-order conditions evaluated at different time instants, and the time inconsistency arises.
Pointed out by \cite{Sethi}, ``a major drawback of the open-loop Stackelberg equilibria is that in general they are not time consistent".
In contrast, according to \cite{Kydland-2,Sethi}, feedback Stackelberg solution is time consistent, i.e., the obtained solution continues to remain optimal at each time instant after the game has begun, and this property is also known as subgame perfect as feedback Stackelberg solution do not depend on system's initial states.

Though the time inconsistency of open-loop Stackelberg solution has been revealed by economist almost half century ago, such type of problems has been seldom investigated within the system control community. 
This paper has the following points to argue the necessity of studying time-consistent open-loop solution of Stackelberg games. Firstly, in the concept of feedback Stackelberg solution, the leader has stagewise advantage over the follower, but not globally; as a constraint, ``such a solution concept requires that the players know the current state of the game in every period" \cite{Sethi}. Yet, open-loop control in deterministic setting is only a function of time and system's initial state, which are clearly known to the players.  Therefore, it is not technically restricted and is natural to announce the leader's open-loop control to the follower before the game begins. Furthermore, for optimal control problems, the set of open-loop controls is the largest admissible control set provided that some constraints are also satisfied to ensure the regularity of controlled systems and cost functionals. Hence, it is very attractive to consider open-loop controls in a Stackelberg game which is indeed composed by two cascaded optimal control problems. This is the case that has been largely investigated in \cite{Basar,Bensoussan,Berkovitz,Freiling,Kydland-2,Shi-1,Simaan-1,Juan-1,Yong-3}.

Secondly,  concerned with the time inconsistency of open-loop Stackelberg solution,
\cite{Kydland-2} points out that ``Faced with this fact, one would expect a great temptation on the part of the dominant player to change his original plan" and ``The feedback solution has the \emph{desirable} characteristic that the plans are intertemporally consistent." Namely, the work \cite{Kydland-2} argues in favour of the feedback solutions as the appropriate solution concept due to its time consistency, where the players should have no rational reason to deviate from the adopted policy in the future. Furthermore, according to the terminology in \cite{Basar}, feedback Stackelberg solution is ``strongly time consistent", as derived by a backward recursion it is regardless at any time instant of previous policies and initial state. Another type of time consistency, called ``weak time consistency", is also introduced in \cite{Basar} to indicate the consistency that the adopted policies will no longer be denied in the future only if past policies are consistent with the equilibrium policies. We are referred to Section 5.6 of \cite{Basar} for more about the two kinds of time consistencies that are compared within the realm of optimal control theory.
As open-loop control depends on system's initial state, the time consistency of open-loop solution of Stackelberg game is likely weaker than that of feedback Stackelberg solution, which merits further investigation.

Thirdly, the time inconsistency of open-loop Stackelberg solution is indeed attributed to the leader's optimal control problem, whose controlled system is a forward-backward difference/differential equation under the follower's best response. In other words, the optimal control problems of forward-backward difference/differential equations are  time-inconsistent, {and yet such kind of time inconsistency has not been investigated before.
%
 %
%
Actually, the reported factors that ruin the time consistency are mainly the non-exponential discounting \cite{Time-discounting, Strotz} and nonlinear terms of conditional expectations \cite{Yong-2} in the objective functionals. Though exponential discounting is of great importance to model people's time preference \cite{Samuelson1937}, empirical researches over the last half century have documented the non-exponential discounting \cite{Time-discounting} that does not possess the property of group or separability any longer.
%
%
Moreover, as there is no nonlinear version of the tower property of conditional expectation, the controller at different time instants is facing different objectives, which are not consistent with the global objective.
In existing literature, there are several different approaches that handle the time inconsistency, and a rule of selecting the preferred solution is called a choice mechanism \cite{Auer}.
The first one is the precommitment choice for which the initial policy is implemented on the lifetime horizon. This approach neglects the time inconsistency, and the optimal policy is optimal only when viewed at the initial time.
Another mechanism is sophisticated/time-consistent choice proposed firstly by Strotz \cite{Strotz}. In the viewpoint of Strotz, the decision maker at different time instants is regarded as different selves, and the time inconsistency suggests a conflict between these different selves.
At any time instant the current self takes account of future selves' decisions, and the equilibrium of this intertemporal game is called a sophisticated policy, or a time-consistent policy.}
%
%
%
%
%
%
Inspired by the work of Strotz, many researchers pay much attention to solving practical problems in economics and finance. Recent years have witnessed the rapid progresses on handling time-inconsistent optimal control problems from the theoretical control community \cite{Hu-1,Hu-2,Ni-1,Ni-2,Yong-2,Wang-1,Yong-1}. 
The study of controlled forward-backward difference/differential equation will extend the boundary of tackling  time-inconsistent optimal control problems. Hence, it is meaningful to investigate time-consistent open-loop solution of Stackelberg games.

In this paper, we propose an open-loop solution concept for a discrete-time deterministic LQ Stackelberg game, which is shown to be weakly time-consistent. {The contents and contributions are listed as follows.}

\begin{itemize}

\item[1.]{A two-tier game framework is introduced to remedy the time inconsistency of open-loop Stackelberg solution. The upper-tier game works according to Stackelberg's scenario, namely, the leader knows the response mapping of the follower but the follower may not know the response mapping of the leader. Two lower-tier intertemporal games are introduced to characterize the follower's and leader's equilibrium response mappings, respectively, which mimic the notion of time-consistent open-loop equilibrium control \cite{Hu-1,Hu-2,Ni-1,Ni-2}.}

    {In other words, the equilibrium of the two-tier game is resorting to solving two unilaterally coupled intertemporal games, by sequentially investigating the time-consistent open-loop equilibrium controls of the cascaded optimal control problems in Stackelberg game. Under the follower's equilibrium response, the leader's controlled system is a forward-backward difference equation (FB$\Delta$E, for short), and the follower's equilibrium response will vary if we just perturb the leader's control action at a single time instant. This deeply distinguishes the second intertemporal game from the one in \cite{Ni-1,Ni-2}, and brings new difficulty to such kind of games.}

\item[2.]{After perturbing the leader's control action, another forward dynamic system is invited to characterize the changed equilibrium response of the follower, besides the one as that in classic variation analysis. Then, an additional adjoint backward equation is introduced to compensate the effect of variation of the follower's equilibrium response. Namely, we need two adjoint backward equations to accomplish the variation analysis of obtaining the stationary condition of the second intertemporal game; see  Proposition \ref{variation} and Theorem \ref{equivalence-1} for more details of this. To the best of the authors, the technique of introducing two adjoint backward equations in  Proposition \ref{variation} has not been seen in existing literature.}

\item[3.]{Then, the equilibrium system of the second intertemporal game includes one forward equation and three backward equations; this differs from the ones for open-loop Stackelberg solution, where the leader's optimal system includes two forward equations and two backward equations. By enlarging the backward state, a modified FB$\Delta$E is introduced. By decoupling this FB$\Delta$E and resorting to the stationary condition, necessary and sufficient conditions are obtained on the unique existence of open-loop equilibrium control of the second intertemporal game. Finally, open-loop equilibrium solution of the two-tier game is derived, which is shown to be weakly time-consistent.}

\end{itemize}

The remainder of the paper is organized as follows. In Section \ref{Section-2-1}, we introduce the Stackelberg game, and investigate its time inconsistency by an approach that is different from the one of \cite{Kydland-2}.  Section \ref{Section-2} introduces the notion of open-loop equilibrium solution of the Stackleberg game, whose full investigations are also presented. Section \ref{Section-4} gives a numerical example, which indicates that the obtained open-loop equilibrium solution is time-consistent. The conclusion is given in Section \ref{Section-5}.

\emph{Notations}. Letting $N$ be a positive integer bigger than 2, denote $\mathbb{T}=\{0,1,\ldots,N-1\}, \mathbb{T}_t=\{t,\ldots,N-1\}$, $\widetilde{\mathbb{T}}_t=\{t,\dots, N\}$ with $t\in \mathbb{T}$. For positive integers $m_1, m_2$, let
\begin{eqnarray}\label{L^2-t-00}
l^2(\mathbb{T}_t; \mathbb{R}^{m})=\Big{\{}{{\{}\rho_k, k\in\mathbb{T}_t{\}}\,\big{|}\, \rho_k\in \mathbb{R}^m, |\rho_k|^2<\infty, k\in\mathbb{T}_t} \Big{\}},~~~m=m_1, m_2.
\end{eqnarray}
{Throughout this paper, a process, say ${\{}a_k, k\in \widetilde{I}{\}}$ valued in some space, may be simply denoted as $a$ when we mention it, where $\widetilde{I}$ is some index set. For example, ${\{}\rho_k, k\in\mathbb{T}_t{\}}$ of (\ref{L^2-t-00}) may be denoted as $\rho$, namely, $\rho\in l^2(\mathbb{T}_t; \mathbb{R}^{m})$.}
If $t$ is replaced by other letters such as $k,\tau$, then $\mathbb{T}_k, \mathbb{T}_\tau, l^2(\mathbb{T}_\tau; \mathbb{R}^{m})$ are similarly defined as $\mathbb{T}_t$ and $l^2(\mathbb{T}_t; \mathbb{R}^{m})$. For any $k\in \mathbb{T}$, let
$$l^2(k; \mathbb{R}^{m})=\{{\rho_k\in \mathbb{R}^{m}\,|\,|\rho_k|^2<\infty}\},~~m=m_1,m_2. $$
For a process, say $v\in l^2(\mathbb{T}_t; \mathbb{R}^{m_2})$, $v|_{\mathbb{T}_k}$ denotes $\{v_k,...,v_{N-1}\}$ with
$k\in \mathbb{T}_t$, {which is the truncation of $v$ on $\mathbb{T}_k$; if $k=t$, $v|_{\mathbb{T}_k}$ will coincide with $v$.  Similar notations such as $\overline{\alpha}^0(x,v)|_{\mathbb{T}_\tau}$ (the truncation of $\overline{\alpha}^0(x,v)$ on $\mathbb{T}_\tau$) appear in the following sections, which are understood from the context.}


\section{Stackelberg game}\label{Section-2-1}

Consider a nonzero-sum deterministic LQ dynamic game associated with the cost functionals
\begin{eqnarray}
\label{cost-1}
J_i(t, x; u, v)=\sum_{k=t}^{N-1}\big{(}X_k^T{Q}_{i}X_k+u_k^T{R}_iu_k+v_k^T{W}_iv_k\big{)}+X_{N}^TG_iX_{N},~~~i=1,2,
%
\end{eqnarray}
which are subject to
\begin{eqnarray}\label{system-dynamics}
\left\{
\begin{array}{l}
X_{k+1}=AX_k+B_1u_k+B_2v_k,\\
X_t=x,~~~k\in \mathbb{T}_t, ~t\in \mathbb{T}.
\end{array}
\right.
\end{eqnarray}
In (\ref{cost-1}) (\ref{system-dynamics}), $\{X_k\in {\mathbb{R}}^n, k\in \widetilde{\mathbb{T}}_t\}$ is the state process, and $u=\{u_k\in \mathbb{R}^{m_1}, k\in \mathbb{T}_t\}, v=\{v_k\in \mathbb{R}^{m_2}, k\in \mathbb{T}_t\}$ are two players' control processes, which are valued in $l^2(\mathbb{T}_t; \mathbb{R}^{m_1})$ and $l^2(\mathbb{T}_t; \mathbb{R}^{m_2})$, respectively.
All the matrices appearing in (\ref{cost-1}) (\ref{system-dynamics}) are deterministic, and the weighting matrices $Q_i, G_i$, $R_i$, $W_i$, $i=1, 2$ are nonnegative definite. 
Corresponding to control processes $u$ and $v$, the players are denoted as Player $u$ and Player $v$, which are minimizing the cost functionals $J_1(t, x; u, v)$ and $J_2(t, x; u, v)$, respectively.
%
%


Stackelberg game, also known as leader-follower game, associated with (\ref{cost-1}) (\ref{system-dynamics}) is formulated in the following way. The leader, Player $v$, first announces her action at the beginning of the game, and  Player $u$ (the follower) seeks his best response strategy so as to minimize $J_1(t,x; u, v)$. Then, anticipating the follower's best response, Player $u$ will select her strategy to minimize $J_2(t,x; u, v)$, and the Stackelberg solution is obtained afterwards. Mathematically, Player $u$ wants to look for a map $\overline{\alpha}^t: \mathbb{R}^n\times l^2(\mathbb{T}_t; \mathbb{R}^{m_2})\mapsto l^2(\mathbb{T}_t; \mathbb{R}^{m_1})$ and then Player $v$ seeks $v^s$ such that
\begin{eqnarray}
%
&&\label{oss-1}
J_1(t, x; \overline{\alpha}^t(x, v), v)=\inf_{u\in l^2(\mathbb{T}_t; \mathbb{R}^{m_1})} J_1(t, x; u, v),\\
&&\label{oss-2}
J_2(t, x;  \overline{\alpha}^t(x, v^s), v^s)=\inf_{v\in l^2(\mathbb{T}_t; \mathbb{R}^{m_2})} J_2(t, x;  \overline{\alpha}^t(x, v), v).
\end{eqnarray}
Noting that (\ref{oss-1}) (\ref{oss-2}) are two cascaded optimal control problems, this exhibits an open-loop Stackelberg solution $(\overline{\alpha}^t(x, v^s), v^s)$, where the leader forces the follower to play in her favor.

The work \cite{Juan-1} studies above Stackelberg game and presents necessary and sufficient condition on the existence and uniqueness of open-loop Stackelberg solution.
For the initial pair $(0, x)$ ($t=0$) and given $v$, problem (\ref{oss-1}) is a standard optimal control problem,
%
%
which is time-consistent, namely, for any $\tau\in \mathbb{T}$, it holds
\begin{eqnarray}\label{Optimal control-follower-2}
J_1(\tau, X_\tau; \overline{\alpha}^0(x,v)|_{\mathbb{T}_\tau}, v|_{\mathbb{T}_\tau})\leq J_1(\tau, X_\tau; u, v|_{\mathbb{T}_\tau}),~~~u\in l^2(\mathbb{T}_\tau, \mathbb{R}^{m_1})
\end{eqnarray}
with $X_\tau$ computed via
\begin{eqnarray*}
\left\{
\begin{array}{l}
X_{k+1}=AX_k+B_1[\overline{\alpha}^0(x,v)]_k+B_2v_k,\\
X_0=x,~~~k\in \mathbb{T}.
\end{array}
\right.
\end{eqnarray*}
Here, $\overline{\alpha}^0(x,v)=([\overline{\alpha}^0(x,v)]_0,...,[\overline{\alpha}^0(x,v)]_{N-1})$ is the unique optimal control that corresponds to the initial pair $(0,x)$ and $v$.
Furthermore, it is shown in \cite{Juan-1} that
\begin{eqnarray}\label{Optimal control-follower-3}
[\overline{\alpha}^0(x,v)]_k=\Phi_{1,k}X_k+\Phi_{2,k}v_k+\Phi_{3,k}\zeta_k,~~~k\in \mathbb{T}
\end{eqnarray}
for some matrices $\Phi_{i,k}, i=1,2,3$; here, $X_k, \zeta_k$ are computed via a FB$\Delta$E, namely, the optimal closed-loop system
\begin{eqnarray}\label{FBDE-1}
\left\{
\begin{array}{l}
X_{k+1}=\bar{A}_kX_k+\bar{B}_{2,k}v_k+\bar{C}_k\zeta_k,\\[1mm]
\zeta_{k}=\bar{D}_k\zeta_{k+1}+\bar{E}_kv_{k+1},\\[1mm]
X_0=x,~~\zeta_N=0,~~k\in \mathbb{T}
\end{array}
\right.
\end{eqnarray}
{for some matrices $\bar{A}_k, \bar{B}_{2,k}, \bar{C}_k, \bar{D}_k, \bar{E}_k, k\in \mathbb{T}$ with compatible dimensions.}
Computing $\zeta_k$ from the backward difference equation, (\ref{FBDE-1}) is rewritten as
\begin{eqnarray}\label{Optimal control-leader}
\left\{
\begin{array}{l}
X_{k+1}=\bar{A}_kX_k+\sum_{\ell=k}^{N-1}\bar{F}_{k,\ell}v_\ell,\\[1mm]
%
%
X_0=x,~~k\in \mathbb{T}.
\end{array}
\right.
\end{eqnarray}
%
%
%
%
%

Under the follower's best response (\ref{Optimal control-follower-3}), the controlled system of problem  (\ref{oss-2})  is the FB$\Delta$E (\ref{FBDE-1}) or equivalently (\ref{Optimal control-leader}).
Substituting $v$ with the optimal one $\widehat{v}^{(0,x)}$, we denote the corresponding optimal state of (\ref{FBDE-1}) by $(\widehat{X}^{(0,x)},\widehat{\zeta}^{(0,x)})$ with $\widehat{X}^{(0,x)}$ the  equilibrium state of this Stackelberg game.
{Hereafter, the superscript $(0,x)$ is added in order to indicate that the concerned variables are corresponding to the initial pair $(0,x)$, which means to differ from those for the time pair $(\tau, y)$ below.}
On the optimal control $\widehat{v}^{(0,x)}$ and by discrete-time variation analysis, the adjoint equation of (\ref{FBDE-1}) is a backward-forward difference equation with state $(\widehat{\alpha}^{(0,x)}, \widehat{\beta}^{(0,x)})$. 
In particular,
\begin{eqnarray}\label{beta-1}
\label{costate-1}
\left\{
\begin{array}{l}
\widehat{\beta}^{(0,x)}_{k+1}=U_{1,k}\widehat{\beta}^{(0,x)}_k+U_{2,k}\widehat{\alpha}^{(0,x)}_k+U_{3,k} \widehat{X}^{(0,x)}_k +U_{4,k}\widehat{v}^{(0,x)}_k+U_{5,k}\widehat{\zeta}^{(0,x)}_k, \\[1mm]
\widehat{\beta}_0^{(0,x)}=0,~~~k\in \mathbb{T},
\end{array}
\right.
\end{eqnarray}
where $U_{p,k}$, $p=1, 2, \ldots, 5$, $k\in \mathbb{T}$, are some deterministic matrices of compatible dimensions.
Also, the superscript $(0,x)$ is indicating that the concerned variables are for the initial pair $(0,x)$. 
%
%
%
%
%
%
%
Furthermore, the unique open-loop Stackelberg solution $(\widehat{u}^{(0,x)}, \widehat{v}^{(0,x)})$ \cite{Juan-1} is
\begin{eqnarray}
\label{follower's optimal control-1}
\left\{
\begin{array}{l}
\widehat{u}^{(0,x)}_k=H_{k}^u\widehat{\xi}^{(0,x)}_k, \\[1mm]
\widehat{v}_k^{(0,x)}=H_{k}^v\widehat{\xi}^{(0,x)}_k,~~~k\in \mathbb{T};
\end{array}
\right.
\end{eqnarray}
{here, $H_{k}^u$, $H_{k}^v, k\in \mathbb{T}$, are functions of $A, B_1, B_2, Q_i, R_i, W_i, G_i (i=1,2)$ with compatible dimensions, and $\widehat{\xi}^{(0,x)}_k={(}(\widehat{\beta}^{(0,x)}_k)^T, (\widehat{X}_k)^T{)}^T$ evolves according to}
\begin{eqnarray}\label{closed-loop-system-0}
\left\{
\begin{array}{l}
\widehat{\xi}^{(0,x)}_{k+1}=S_k\widehat{\xi}^{(0,x)}_k,\\[1mm]
\widehat{\xi}^{(0,x)}_0=(0^T, x^T)^T,~~k\in \mathbb{T}
\end{array}
\right.
\end{eqnarray}
for some matrices $S_k, k\in \mathbb{T}$. {It is worth mentioning that the initial value $\widehat{\beta}^{(0,x)}_0$ is 0.}

At $\tau\in \mathbb{T}$ ($\tau>0$), the equilibrium state is $\widehat{X}^{(0,x)}_\tau$, {and to this end it is denoted as $y$ for notation simplicity. Reconsider this Stackelberg game at the initial pair $(\tau, y)$, and the corresponding variables are labelled by the superscript $(\tau, y)$ to distinguish those for $(0,x)$.
In this case,} the controlled system of the leader is
\begin{eqnarray}\label{FBDE-tau}
\left\{
\begin{array}{l}
X_{k+1}=\bar{A}_kX_k+\bar{B}_{2,k}v_k+\bar{C}_k\zeta_k,\\[1mm]
\zeta_{k}=\bar{D}_k\zeta_{k+1}+\bar{E}_kv_{k+1},\\[1mm]
X_\tau=y,~~\zeta_N=0,~~k\in \mathbb{T}_\tau,
\end{array}
\right.
\end{eqnarray}
and the adjoint equation of optimal backward equation is
\begin{eqnarray}
\label{costate-tau}
\left\{
\begin{array}{l}
\widehat{\beta}^{(\tau,y)}_{k+1}=U_{1,k+1}\widehat{\beta}^{(\tau,y)}_k+U_{2,k+1}\widehat{\alpha}^{(\tau,y)}_k+U_{3,k+1}\widehat{X}^{(\tau,y)}_k +U_{4,k+1}\widehat{v}^{(\tau,y)}_k+U_{5,k+1}\widehat{\zeta}^{(\tau,y)}_k, \\[1mm]
\widehat{\beta}_\tau^{(\tau,y)}=0,~~~k\in \mathbb{T}.
\end{array}
\right.
\end{eqnarray}
Furthermore, the leader's equilibrium control (corresponding to the initial pair $(\tau, y)$) is
\begin{eqnarray*}
\widehat{v}^{(\tau,y)}_k=H_{k}^v\widehat{\xi}^{(\tau,y)}_k,~~~k\in \mathbb{T}_\tau
\end{eqnarray*}
with $\widehat{\xi}^{(\tau,y)}_k={(}(\widehat{\beta}^{(\tau,y)}_k)^T, (\widehat{X}^{(\tau,y)}_k)^T{)}^T$, which evolves according to
\begin{eqnarray}
\left\{
\begin{array}{l}
\xi^{(\tau,y)}_{k+1}=S_k\xi^{(\tau,y)}_k,\\[1mm]
\xi^{(\tau,y)}_\tau=(0^T, y^T)^T,~~~k\in \mathbb{T}_\tau.
\end{array}
\right.
\end{eqnarray}
Note that the initial value $\beta^{(\tau,y)}_\tau$ is 0. On the other hand, following (\ref{closed-loop-system-0}), $\beta^{(0,x)}_\tau$ is generally nonzero provided that $x\neq 0$, i.e., $\beta^{(\tau,y)}_\tau\neq \beta^{(0,x)}_\tau.$
Hence, $\xi^{(\tau,y)}_k\neq \xi^{(0,x)}_k, k\in \mathbb{T}_\tau$, which implies
\begin{eqnarray}\label{TI-1}
\widehat{v}^{(\tau,y)}_k\neq \widehat{v}^{(0,x)}_k,~~~  k\in \mathbb{T}_\tau.
\end{eqnarray}
%
%
%
%
%
Namely, the truncation of open-loop Stackelberg solution (for the initial pair $(0, x)$) on  $\mathbb{T}_\tau$ is not the open-loop Stackelberg solution for the initial pair $(\tau, \widehat{X}^{(0,x)}_\tau)$. This phenomenon is termed the time inconsistency of open-loop Stackelberg solution.
Essentially, above derivation from (\ref{Optimal control-follower-2}) to (\ref{TI-1}) shows that optimal control problem (\ref{oss-2}) associated with the the FB$\Delta$E (\ref{FBDE-1}) is time-inconsistent!
The following is a simple numerical example.

\begin{example}
In (\ref{cost-1}) (\ref{system-dynamics}), let $A= B_1=B_2=Q_1=R_{1}=W_{2}=G_{1}=1$, $R_{2}=W_{1}=0$, $Q_2=3$, $G_{2}=2$, $t=0$, $N=3$, and $x=1$. Check the time inconsistency of open-loop Stackelberg solution.

\end{example}

\emph{Solution:} Simple calculations give the unique open-loop Stackelberg solution for the initial pair $(0, x)$:
\begin{eqnarray*}
&&(\widehat{u}^{(0,x)}_0, \widehat{v}^{(0,x)}_0)=(-0.4240, -0.3363),~~~(\widehat{u}^{(0,x)}_1, \widehat{v}^{(0,x)}_1)=(-0.1843, 0.0465),\\
&&(\widehat{u}^{(0,x)}_2, \widehat{v}^{(0,x)}_2)=(-0.0823, 0.0626).
\end{eqnarray*}
At time instant 1, the equilibrium state $\widehat{X}^{(0,x)}_1=0.2397$, which is denoted as $y$. Now, reconsider this Stackelberg game for the initial pair $(1, y)$ whose open-loop Stackelberg solution is
\begin{eqnarray*}
(\widehat{u}_1^{(1,y)}, \widehat{v}^{(1,y)}_1)=(-0.0942, -0.0856), ~~(\widehat{u}^{(1,y)}_2, \widehat{v}^{(1,y)}_2)=(-0.0342, 0.0086).
\end{eqnarray*}
Clearly,
\begin{eqnarray*}
(\widehat{u}^{(1,y)}_1, \widehat{v}^{(1,y)}_1)\neq(\widehat{u}^{(0,x)}_1, \widehat{v}^{(0,x)}_1), ~~(\widehat{u}_2^{(1,y)}, \widehat{v}^{(1,y)}_2)\neq(\widehat{u}_2^{(0,x)}, \widehat{v}^{(0,x)}_2),
 \end{eqnarray*}
and the open-loop Stackelberg solution is time-inconsistent.   \hfill $\square$

To conclude this section, finding open-loop Stackelberg solution is indeed divided into two cascaded optimal control problems (\ref{oss-1}) (\ref{oss-2}), 
%
%
and it has been shown that under the follower's best response the leader's optimal control problem (\ref{oss-2}) is time-inconsistent.
%
%
%
Hence, the overall cascaded optimal control problems (\ref{oss-1}) (\ref{oss-2}) are time-inconsistent, although for a given $v$ problem (\ref{oss-1}) is time-consistent. 
%
%
%
 %
So, if the equilibrium solution $(\widehat{u}^{(0, x)}, \widehat{v}^{(0,x)})$ is selected once for all, then the story ends without subsequent replanning; {here, $(\widehat{u}^{(0, x)}, \widehat{v}^{(0,x)})$ is named the  precommitted solution \cite{Auer}}. If instead, policy choice is sequent and made period after period, then the leader has an incentive to deviate from his initial policy $\widehat{v}^{(0,x)}$ later on as it is no longer optimal! In this case, proper notion on Stackelberg equilibrium needs to be deliberately selected in order to remedy the time inconsistency.

\section{Open-loop equilibrium solution}\label{Section-2}

\subsection{Definition}

{In this paper, we adopt Strotz's philosophy to handle the time inconsistency of Stackelberg game. Let us firstly say more on the one-player time-inconsistent optimal control. Facing the time inconsistency, a sophisticated decisionmaker is unable to precommit to the strategy selected at present, and his future selves may deviate from what he makes now.
Aware of this conflict, \cite{Pollak,Strotz} propose the concept of consistent planning, where the decisionmaker anticipates his strategy in the future and act today accordingly; more specifically, a time-consistent solution is introduced as the equilibrium outcome of an intertemporal game between different selves (individuals) who make decisions sequentially over time.
%
%
For the time-consistent open-loop equilibrium control, we are referred to, for instance \cite{Ni-1}, where a time-inconsistent mean-field stochastic LQ problem is studied.
In a word, Strotz's philosophy resorts to a (one-tier) intertemporal game to resolve the time inconsistency of the optimal control problem. Here, the term ``one-tier" is deliberately used in order to distinguish the one ``two-tier" below for Stackelberg game.}

{We again consider a Stackelberg scenario: Player $v$ is the leader and Player $u$ the follower. To remedy the time inconsistency of (\ref{oss-1}) (\ref{oss-2}), we instead
consider an open-loop control pair for which no deviation at a single time
instant let each player incurs a lower cost, that is for the follower
changing his response at a single time instant considering the strategy of the
leader unchanged, and for the leader to change his leading strategy at
a single time instant considering the response of the follower.
%
Namely, we are facing a two-tier game: besides the upper-tier game between the leader and follower, there are two lower-tier intertemporal games among the selves of the leader and of the follower, respectively.}

\begin{definition}\label{definition}

Concerned with the Stackelberg game associated (\ref{cost-1}) (\ref{system-dynamics}), a pair $(u^*,v^*)\in l^2(\mathbb{T}_t; \mathbb{R}^{m_1})\times l^2(\mathbb{T}_t; \mathbb{R}^{m_2})$ provides the unique open-loop equilibrium solution for the initial pair $(t,x)$, if
\begin{itemize}
\item[i)] For each $v\in l^2(\mathbb{T}_t; \mathbb{R}^{m_2})$, there exists a unique $\bar{u}\in l^2(\mathbb{T}_t; \mathbb{R}^{m_1})$ such that
    \begin{eqnarray}\label{inequality-0}
      J_1(k, X_k; \bar{u}|_{\mathbb{T}_k}, v|_{\mathbb{T}_k})\leq J_1(k, X_k; (u_k, \bar{u}|_{\mathbb{T}_{k+1}}), v|_{\mathbb{T}_k})
    \end{eqnarray}
holds for any $k\in\mathbb{T}_t$ and any $u_k\in l^2(k; \mathbb{R}^{m_1})$. Equivalently, there exists a unique map $\alpha^t: \mathbb{R}^n\times l^2(\mathbb{T}_t; \mathbb{R}^{m_2})\mapsto l^2(\mathbb{T}_t; \mathbb{R}^{m_1})$ such that
    \begin{eqnarray}\label{inequality-1}
      J_1(k, X_k; {\alpha^t}(x, v)|_{\mathbb{T}_k}, v|_{\mathbb{T}_k})\leq J_1(k, X_k; (u_k, {\alpha^t}(x, v)|_{\mathbb{T}_{k+1}}), v|_{\mathbb{T}_k})
    \end{eqnarray}
    holds for each $v\in l^2(\mathbb{T}_t; \mathbb{R}^{m_2})$, any $k\in\mathbb{T}_t$ and any $u_k\in l^2(k; \mathbb{R}^{m_1})$ with  
    %
    $\alpha^t(x, v)=([\alpha^t(x, v)]_t,\ldots, [\alpha^t(x, v)]_{N-1})$.
    The state $X_k$ in (\ref{inequality-1}) is computed via
    \begin{eqnarray}\label{system-dynamics-2}
\left\{
\begin{array}{l}
X_{k+1}=AX_k+B_1[\alpha^t(x, v)]_{k}+B_2v_k,\\
X_t=x,~~~~k\in \mathbb{T}_t.
\end{array}
\right.
\end{eqnarray}

\item[ii)]  There exists a unique $v^*\in l^2(\mathbb{T}_t; \mathbb{R}^{m_2})$ such that
    \begin{eqnarray}\label{inequality-2}
    J_2(k, X_k^*; {\alpha^t}(x, v^*)|_{\mathbb{T}_k}, v^*|_{\mathbb{T}_k})\leq J_2(k, X_k^*; {\alpha^t}(x, v^{*-k})|_{\mathbb{T}_k}, (v_k, v^*|_{\mathbb{T}_{k+1}}))
    \end{eqnarray}
holds for any $k\in\mathbb{T}_t$ and any $v_k\in l^2(k; \mathbb{R}^{m_2})$. The state $X_k^*$ in (\ref{inequality-2}) is computed via
\begin{eqnarray}\label{system-dynamics-3}
\left\{
\begin{array}{l}
X^*_{k+1}=AX^*_k+B_1[\alpha^t(x, v^*)]_{k}+B_2v^*_k,\\
X^*_t=x,~~~~k\in \mathbb{T}_t,
\end{array}
\right.
\end{eqnarray}
and the internal state of $J_2(k, X_k^*; {\alpha^t}(x, v^{*-k})|_{\mathbb{T}_k}, (v_k, v^*|_{\mathbb{T}_{k+1}}))$ is
\begin{eqnarray}\label{internal-state}
\left\{
\begin{array}{l}
X_{\ell+1}=AX_\ell+B_1[\alpha^t(x, v^{*-k})]_{\ell}+B_2v^{*-k}_\ell,\\
X_k=X_k^*,~~\ell \in \mathbb{T}_k,
\end{array}
\right.
\end{eqnarray}
where
\begin{eqnarray*}
v^{*-k}_\ell=\left\{
\begin{array}{ll}
v_k,  ~~&\ell=k,\\
v^*_\ell,  ~~&\ell\neq k, \ell \in \mathbb{T}_t.
\end{array}
\right.
\end{eqnarray*}

\item[iii)]  $u^*=\alpha^t(x, v^*)$.
\end{itemize}

\end{definition}

\begin{remark}
{
The three steps in Definition \ref{definition} describe a two-tier game. The upper-tier game works in Stackelberg' scenario between the leader $v$ and follower $u$, namely, the leader knows the response mapping of the follower but the follower may not know the response mapping of the leader. In the lower tier, two intertemporal games, see (\ref{inequality-1}) (\ref{inequality-2}), are introduced to remedy the time inconsistency of open-loop Stackelberg solution. For fixed $v$, $\alpha^t(x, v)$ of (\ref{inequality-1}) is the time-consistent open-loop equilibrium control of problem (\ref{oss-1}), which is shown to have the form (\ref{follower's optimal control}). Consider the other intertemporal game related to (\ref{inequality-2}). If at time instant $k$ we select $v_k$ instead of $v^*_k$, ${\alpha^t}(x, v^{*-k})|_{\mathbb{T}_k}$ will differ from ${\alpha^t}(x, v^{*})|_{\mathbb{T}_k}$. Due to (\ref{follower's optimal control}) below, the state of
\begin{eqnarray}\label{system-leader-0}
\left\{
\begin{array}{l}
X_{\ell+1}=\widetilde{A}_\ell X_\ell+\widetilde{B}_\ell v^{*-k}_{\ell}+\widetilde{C}_\ell \pi_{\ell+1},\\[1mm]
\pi_\ell=C_\ell^Tv^{*-k}_\ell+\widetilde{A}_\ell^T\pi_{\ell+1},\\[1mm]
X_t=x,~~~\pi_{N}=0,~~~\ell\in \mathbb{T}_t
\end{array}
\right.
\end{eqnarray}
is used to compute ${\alpha^t}(x, v^{*-k})$. Therefore, we need three dynamic systems (\ref{system-dynamics-3}) (\ref{internal-state}) (\ref{system-leader-0}) to characterize (\ref{inequality-2}). This makes the concerned problem much complicated, and such an intertemporal game differs significantly from the one for time-inconsistent optimal control \cite{Ni-1, Ni-2}.
}

\end{remark}

\subsection{Characterizations}


%
%
%

%

We firstly characterize of the map $\alpha^t$ of $i)$ of Definition \ref{definition}. The following result gives conditions on the existence and uniqueness of $\alpha^t$, whose proof follows directly from Theorem III.5 of \cite{Ni-1} and is omitted here.

\begin{theorem}
The following statements are equivalent.
\begin{itemize}
\item[i)] There exists a unique map $\alpha^t$ such that (\ref{inequality-1}) holds for each $v\in l^2(\mathbb{T}_t; \mathbb{R}^{m_2})$, any $k\in\mathbb{T}_t$ and any $u_k\in l^2(k; \mathbb{R}^{m_1})$.
\item[ii)] The matrices
\begin{eqnarray}\label{M}
 M_k=B_1^TP_{k+1}B_1+R_{1}, ~~~k\in k\in\mathbb{T}_t
\end{eqnarray}
are all positive definite, where $P_{k+1}$ is computed via
\begin{eqnarray}\label{P}
\left\{
\begin{array}{l}
P_k=Q_1+A^TP_{k+1}A-A^TP_{k+1}B_1M_k^{-1}B_1^TP_{k+1}A,\\[1mm]
P_{N}=G_1,~~~k\in \mathbb{T}_t.
\end{array}
\right.
\end{eqnarray}

\end{itemize}

In this case, the value of $\alpha^t(x, v)$ is given by
\begin{eqnarray}\label{follower's optimal control}
[\alpha^t(x, v)]_k=-\big{(}H^1_kX_k+H^2_kv_k+H_k^3\pi_{k+1}\big{)},~~k\in\mathbb{T}_t,
\end{eqnarray}
where
\begin{eqnarray*}
\left\{
\begin{array}{l}
H_{k}^1=M_k^{-1}B_1^TP_{k+1}A, \\[1mm]
H_k^2=M_k^{-1}B_1^TP_{k+1}B_2, \\[1mm]
H_k^3=M_k^{-1}B_1^T, ~~~k\in\mathbb{T}_{t},
\end{array}
\right.
\end{eqnarray*}
and $X_k$, $\pi_k$, $k\in\mathbb{T}_t$, are computed via
\begin{eqnarray}\label{system-leader}
\left\{
\begin{array}{l}
X_{k+1}=\widetilde{A}_kX_k+\widetilde{B}_kv_{k}+\widetilde{C}_k\pi_{k+1},\\[1mm]
\pi_k=C_k^Tv_k+\widetilde{A}_k^T\pi_{k+1},\\[1mm]
X_t=x,~~~\pi_{N}=0,~~~k\in \mathbb{T}_t
\end{array}
\right.
\end{eqnarray}
with
\begin{eqnarray*}
\left\{
\begin{array}{l}
\widetilde{A}_k=A-B_1H_k^1,\\[1mm]
\widetilde{B}_k=B_2-B_1H_k^2,\\[1mm]
\widetilde{C}_k=-B_1H_k^3,\\[1mm]
C_k=(B_2^T-B_2^TP_{k+1}B_1M_k^{-1}B_1^T)P_{k+1}A,\\[1mm]
k\in\mathbb{T}_t.
\end{array}
\right.
\end{eqnarray*}

\end{theorem}

Under (\ref{follower's optimal control}), we now characterize the control process $v^*$ of $ii)$ of Definition \ref{definition}. In this case, the controlled system of the leader is (\ref{system-leader}),
which is a FB$\Delta$E. The following proposition gives the expression of perturbed cost functional, whose proof together with those of  Theorem \ref{equivalence-1} and  Theorem \ref{Theorem-stationary-condition} are given in Section \ref{proof}.

\begin{proposition}
\label{variation}

For $v \in l^2(\mathbb{T}_t; \mathbb{R}^{m_2})$, $\widetilde{v}_k \in l^2(k; \mathbb{R}^{m_2})$ and $\varepsilon \in\mathbb{R}$, let
\begin{eqnarray*}
v^\varepsilon_\ell=\left\{
\begin{array}{ll}
v_k+\varepsilon \widetilde{v}_k,&~~~\ell=k,\\
v_\ell,&~~~\ell\neq k, \ell\in \mathbb{T}_t,
\end{array}
\right.
\end{eqnarray*}
and
$v^\varepsilon=(v_t^\varepsilon, \ldots, v^\varepsilon_{N-1})$.
Then, it holds that
\begin{eqnarray}
\label{difference-formula-5}
&&J_2(k, X_k; \alpha^t(x, v^\varepsilon)|_{\mathbb{T}_{k}}, (v_k+\varepsilon \widetilde{v}_k,  v|_{\mathbb{T}_{k+1}}))-
J_2(k, X_k; \alpha^t(x, v)|_{\mathbb{T}_{k}}, v|_{\mathbb{T}_{k}})\nonumber\\
&&=2\varepsilon \Big{[}\widetilde{B}_k^TZ_{k+1}+\widetilde{B}_{k}^T\overline{Z}_{k+1}^{(k)}+W_2v_k
+(H_k^2)^TR_2(H_k^1X_k+H_k^2v_k+H_k^3\pi_{k+1}) \nonumber\\
&&\hphantom{=}+\sum_{i=t}^{k-1}C_{k}\widetilde{A}_{k-1}\cdots \widetilde{A}_{i+1}\widetilde{C}_{i}^T\overline{Z}_{i+1}^{(k)}
\Big{]}^T\bar{v}_k+\varepsilon^2\widehat{J}_2(k, 0; \widetilde{v}_k).
\end{eqnarray}
Here,
$Z_{k+1}$, $\overline{Z}_{k+1}^{(k)}, \cdots, \overline{Z}_{t+1}^{(k)}$  are  computed via
\begin{eqnarray}
\left\{
\begin{array}{l}
Z_k=Q_2X_{k}+A^TZ_{k+1}, \\[1mm]
Z_{N}=G_2X_N,~~k \in \mathbb{T}_t,
\end{array}
\right.
\end{eqnarray}
and
\begin{eqnarray}
\label{backwards equation-2}
\left\{
\begin{array}{l}
\left\{
\begin{array}{l}
\overline{Z}^{(k)}_{\ell}=(H_{\ell}^1)^TR_2(H_{\ell}^1X_{\ell}+H^2_{\ell}v_{\ell}+H^3_{\ell}
\pi_{\ell+1})-(H_{\ell}^1)^TB_1^TZ_{\ell+1}+\widetilde{A}_{\ell}^T\overline{Z}^{(k)}_{\ell+1}, \\[1mm]
\overline{Z}^{(k)}_{i}=\widetilde{A}_{i}^T\overline{Z}^{(k)}_{i+1}, \\[1mm]
\overline{Z}^{(k)}_{N}=0, \\[1mm]
\ell\in\mathbb{T}_{k},~~i\in {\{}t, t+1, \ldots,k-1{\}},\\[1mm]
\end{array}
\right.\\
k \in \mathbb{T}_t
\end{array}
\right.
\end{eqnarray}
with
\begin{eqnarray}
\label{backwards equation-21}
\left\{
\begin{array}{l}
X_{k+1}=\widetilde{A}_kX_k+\widetilde{B}_kv_{k}+\widetilde{C}_k\pi_{k+1}, \\[1mm]
\pi_k=C_k^Tv_k+\widetilde{A}_k^T\pi_{k+1}, \\[1mm]
X_{t}=x,~~\pi_{N}=0, \\[1mm]
k\in\mathbb{T}_{t}.
\end{array}
\right.
\end{eqnarray}
Furthermore, $\widehat{J}_2(k, 0; \widetilde{v}_k)$ of (\ref{difference-formula-5}) is given by
\begin{eqnarray}
\label{convex_2}
&&\widehat{J}_2(k, 0; \widetilde{v}_k)=\sum_{\ell=k}^{N-1}\xi_\ell^TQ_2\xi_\ell+
\sum_{\ell=k+1}^{N-1}(\eta^{(k)}_{\ell})^T(H_{\ell}^1)^TR_2H_{\ell}^1\eta^{(k)}_{\ell}
+ \widetilde{v}_k^TW_{2} \widetilde{v}_k \nonumber\\
&&\hphantom{\widehat{J}_2(k, 0; \bar{v}_k)=}
+\xi_{N}^TG_2\xi_{N}+(H_k^1\eta^{(k)}_k+H^2_k \widetilde{v}_k)^TR_2(H_k^1\eta^{(k)}_k+H^2_k \widetilde{v}_k)
\end{eqnarray}
with
\begin{eqnarray}
\label{Y-1}
\left\{
\begin{array}{l}
\eta^{(k)}_{\ell+1}=\widetilde{A}_\ell\eta^{(k)}_{\ell}, \\[1mm]
\eta^{(k)}_{k+1}=\widetilde{A}_k\eta^{(k)}_k+ \widetilde{B}_k\widetilde{v}_k, \\[1mm]
\eta^{(k)}_{i+1}=\widetilde{A}_{i}\eta^{(k)}_{i}+ \widetilde{C}_{i} \widetilde{A}_{i+1}^T\cdots \widetilde{A}_{k-1}^T%
C_{k}^T\widetilde{v}_k,\\[1mm]
\eta^{(k)}_{t}=0,\\[1mm]
\ell\in \mathbb{T}_{k+1},~~~i\in \{t, t+1, \ldots, k-1\},
\end{array}
\right.
\end{eqnarray}
and
\begin{eqnarray}
\label{Y-2}
\left\{
\begin{array}{l}
\xi_{\ell+1}=A\xi_\ell-B_1H_{\ell}^1\eta^{(k)}_{\ell}, \\[1mm]
\xi_{k+1}=A\xi_k-B_1H_k^1\eta^{(k)}_k+ \widetilde{B}_k\widetilde{v}_k, \\[1mm]
\xi_k=0, ~~~\ell\in \mathbb{T}_{k+1}.
\end{array}
\right.
\end{eqnarray}
In (\ref{Y-1}), the following notation is adopted
\begin{eqnarray*}
\widetilde{A}_{i+1}^T\cdots \widetilde{A}_{k-1}^T=\left\{
\begin{array}{ll}
I,&i+1>k-1,\\[1mm]
\widetilde{A}_{k-1}^T,&i+1=k-1,\\[1mm]
\widetilde{A}_{i+1}^T\cdots \widetilde{A}_{k-1}^T,&i+1<k-1,
\end{array}
\right.
\end{eqnarray*}
and $\widetilde{A}_{k-1}\cdots \widetilde{A}_{i+1}$ of (\ref{difference-formula-5}) is the {transpose} of $\widetilde{A}_{i+1}^T\cdots \widetilde{A}_{k-1}^T$.

\end{proposition}

\begin{remark}

Due to the expression (\ref{follower's optimal control}), $[\alpha^t(x,v^\varepsilon)]_t,...,[\alpha^t(x,v^\varepsilon)]_{N-1}$ are all modified if we just replace $v_k$ by $v_k^\varepsilon$. Note that $\eta^{(k)}$ is introduced to characterize the difference between $X$ and $X^\varepsilon$ that are given in (\ref{system-leader}) and
\begin{eqnarray}\label{system-leader-e}
\left\{
\begin{array}{l}
X^\varepsilon_{k+1}=\widetilde{A}_kX^\varepsilon_k+\widetilde{B}_kv^\varepsilon_{k}
+\widetilde{C}_k\pi^\varepsilon_{k+1}, \\[1mm]
\pi^\varepsilon_k=C_k^Tv^\varepsilon_k+\widetilde{A}_k^T\pi^\varepsilon_{k+1}, \\[1mm]
X^\varepsilon_t=x,~~~\pi^\varepsilon_{N}=0, \\[1mm]
k\in\mathbb{T}_{t}.
\end{array}
\right.
\end{eqnarray}
Furthermore, the backward equation on $\overline{Z}^{(k)}$ is introduced to cancel the effect of $\eta^{(k)}$ (in (\ref{difference formula-3})). To the best of the authors, the technique of introducing two adjoint equations in Proposition \ref{variation} has not been seen in existing literature.

\end{remark}

Through (\ref{backwards equation-2}), one gets
\begin{eqnarray*}
&&\sum_{i=t}^{k-1}C_{k}\widetilde{A}_{k-1}\cdots\widetilde{A}_{i+1}\widetilde{C}_{i}^T\overline{Z}_{i+1}^{(k)} \nonumber\\
&&=\sum_{i=t}^{k-1}C_{k}\widetilde{A}_{k-1}\cdots\widetilde{A}_{i+1}\widetilde{C}_{i}^T\widetilde{A}_{i+1}^T
\cdots\widetilde{A}_{k-1}^T\overline{Z}_{k}^{(k)} \nonumber\\
&&=\sum_{i=t}^{k-1}C_{k}\widetilde{A}_{k-1}\cdots\widetilde{A}_{i+1}\widetilde{C}_{i}^T\widetilde{A}_{i+1}^T
\cdots\widetilde{A}_{k-1}^T\Big{[}(H_{k}^1)^TR_2(H_{k}^1X_{k}+H^2_{k}v_{k}+H_k^3
\pi_{k+1}) \nonumber\\
&&\hphantom{=}-(H_{k}^1)^TB_1^TZ_{k+1}+\widetilde{A}_{k}^T\overline{Z}_{k+1}^{(k)}\Big{]} \nonumber\\
&&=\sum_{i=t}^{k-1}D_i^{(k)}\Big{[}(H_{k}^1)^TR_2(H_{k}^1X_{k}+H^2_{k}v_{k}+H_k^3
\pi_{k+1})-(H_{k}^1)^TB_1^TZ_{k+1}+\widetilde{A}_{k}^T\overline{Z}_{k+1}^{(k)}\Big{]},
\end{eqnarray*}
where $D_i^{(k)}=C_k\widetilde{A}_{k-1}\cdots \widetilde{A}_{i+1}\widetilde{C}_i^T\widetilde{A}_{i+1}^T \cdots \widetilde{A}_{k-1}^T$,~~$i \in {\{}t, t+1, \ldots k-1 {\}}$.
Then, the following result holds

\begin{lemma}
\label{stationary condition}

The equalities
\begin{eqnarray}
&&\widetilde{B}_k^TZ_{k+1}+\widetilde{B}_{k}^T\overline{Z}_{k+1}^{(k)}+W_2v_k
+(H^2_k)^TR_2(H_k^1X_k+H^2_kv_k+H_k^3\pi_{k+1}) \nonumber\\
&&+\sum_{i=t}^{k-1}C_{k}\widetilde{A}_{k-1}\cdots \widetilde{A}_{i+1}\widetilde{C}_{i}^T\overline{Z}_{i+1}^{(k)} \nonumber\\
&&=\Big{[}(H^2_k)^TR_2H_k^1
+\sum_{i=t}^{k-1}D_i^{(k)}(H_k^1)^TR_2H_k^1\Big{]}X_k
+\Big{[}W_2+(H^2_k)^TR_2H^2_k+\sum_{i=t}^{k-1}D_i^{(k)}(H_k^1)^TR_2H^2_k\Big{]}v_k \nonumber\\
&&\hphantom{=}+\Big{[}\widetilde{B}_{k}^T-\sum_{i=t}^{k-1}D_i^{(k)}(H_k^1)^TB_1^T\Big{]}Z_{k+1}+(\widetilde{B}_{k}^T
+\sum_{i=t}^{k-1}D_i^{(k)}\widetilde{A}_k^T)\overline{Z}_{k+1}^{(k)} \nonumber\\
&&\hphantom{=}+\Big{[}(H^2_k)^TR_2H_k^3+\sum_{i=t}^{k-1}D_i^{(k)}(H_k^1)^TR_2H_k^3\Big{]}\pi_{k+1},~~~~~k\in \mathbb{T}_t
\end{eqnarray}
are satisfied.

\end{lemma}

{By Proposition \ref{variation} and Lemma \ref{stationary condition}, we can characterize the stationary condition of optimization problem related to (\ref{inequality-2}) that is indeed the first-order optimality condition. As the weighting matrices $Q_i, G_i, R_i, W_i, i=1,2$, are all nonnegative definite, the considered optimization problem is convex. Hence, the stationary condition is necessary and sufficient to characterize the minimizer of (\ref{inequality-2}).
}

\begin{theorem}
\label{equivalence-1}

For the initial pair $(t,x)$, the following statements are equivalent.
\begin{itemize}
\item[i)]  There exists a unique $v^*\in l^2(\mathbb{T}_t; \mathbb{R}^{m_2})$ such that (\ref{inequality-2}) holds for any $k\in\mathbb{T}_t$ and any $v_k\in l^2(k; \mathbb{R}^{m_2})$.
\item[ii)]  There exists a unique $v^*\in l^2(\mathbb{T}_t; \mathbb{R}^{m_2})$ such that the stationary condition
\begin{eqnarray}
\label{stationary-condition-1}
&&0=\Big{[}\widetilde{B}_{k}^T-\sum_{i=t}^{k-1}D_i^{(k)}(H_k^1)^TB_1^T\Big{]}Z_{k+1}^*+(\widetilde{B}_{k}^T
+\sum_{i=t}^{k-1}D_i^{(k)}\widetilde{A}_k^T)\overline{Z}_{k+1}^{(k)*}+\Big{[}(H^2_k)^TR_2H_k^1\nonumber\\
&&\hphantom{0=}+\sum_{i=t}^{k-1}D_i^{(k)}(H_k^1)^TR_2H_k^1\Big{]}X_k^* +\Big{[}W_2+(H^2_k)^TR_2H^2_k+\sum_{i=t}^{k-1}D_i^{(k)}(H_k^1)^TR_2H^2_k\Big{]}v_k^*\nonumber\\
&&\hphantom{0=}
+\Big{[}(H^2_k)^TR_2H_k^3+\sum_{i=t}^{k-1}D_i^{(k)}(H_k^1)^TR_2H_k^3\Big{]}\pi_{k+1}^*,~~~k\in \mathbb{T}_t
\end{eqnarray}
holds. Here, 
$Z_{k+1}^*$ and $\overline{Z}_{k+1}^{(k)*}$ are computed via
\begin{eqnarray}
\label{backwards equation-21}
\left\{
\begin{array}{l}
Z_k^*=Q_2X_{k}^*+A^TZ_{k+1}^*, \\[1mm]
Z_{N}^*=G_2X_N^*, ~~~~k\in\mathbb{T}_{t},
\end{array}
\right.
\end{eqnarray}
and
\begin{eqnarray}
\label{backwards equation-2-21}
\left\{
\begin{array}{l}
\left\{
\begin{array}{l}
\overline{Z}_{\ell}^{(k)*}=(H_{\ell}^1)^TR_2(H_{\ell}^1X_{\ell}^*+H^2_{\ell}v_{\ell}^*+H_{\ell}^3
\pi_{\ell+1}^*)-(H_{\ell}^1)^TB_1^TZ_{\ell+1}^*+\widetilde{A}_{\ell}^T\overline{Z}_{\ell+1}^{(k)*}, \\[1mm]
\overline{Z}_{N}^{(k)*}=0, ~~~~\ell\in\mathbb{T}_{k}, \\[1mm]
\end{array}
\right. \\
k \in \mathbb{T}_t
\end{array}
\right.
\end{eqnarray}
with
\begin{eqnarray}
\label{state-X*}
\left\{
\begin{array}{l}
X_{k+1}^*=\widetilde{A}_kX_k^*+\widetilde{B}_kv_{k}^*+\widetilde{C}_k\pi_{k+1}^*, \\[1mm]
\pi_k^*=C_k^Tv_k^*+\widetilde{A}_k^T\pi_{k+1}^*, \\[1mm]
X_t^*=x,~~\pi_{N}^*=0, \\[1mm]
k\in\mathbb{T}_{t}.
\end{array}
\right.
\end{eqnarray}

In this case, $v^*$ of ii) is the one of i).

\end{itemize}

\end{theorem}

\begin{remark}
From (\ref{backwards equation-2-21}), one has
\begin{eqnarray*}
\overline{Z}^{(k_1)*}_{\ell}=\overline{Z}^{(k_2)*}_{\ell},~~~\forall k_1, k_2\in \mathbb{T}_t,~k_1<k_2,~\ell\in \mathbb{T}_{k_2}.
\end{eqnarray*}
Hence, (\ref{backwards equation-2-21}) is simplified to
\begin{eqnarray}
\label{backwards equation-2-211}
\left\{
\begin{array}{l}
\overline{Z}_{k}^{*}=(H_{k}^1)^TR_2(H_{k}^1X_{k}^*+H^2_{k}v_{k}^*+H^3_k
\pi_{k+1}^*)-(H_{k}^1)^TB_1^TZ_{k+1}^*+\widetilde{A}_{k}^T\overline{Z}_{k+1}^{*}, \\[1mm]
\overline{Z}_{N}^{*}=0, ~~~~k\in\mathbb{T}_{t}.
\end{array}
\right.
\end{eqnarray}
\end{remark}

%
%

%
%
Introduce the following notations
\begin{eqnarray}
\label{dimension-rising-matrix}
\left\{
\begin{array}{l}
\mathbf{Z}_k^{*}=\left(\begin{array}{cc} Z_k^* \\[1mm] \overline{Z}_k^{*} \\ \pi_k^* \end{array}\right),
\mathbf{H}_{k}=\left(\begin{array}{cc} Q_2 \\[1mm] (H_{k}^1)^TR_2H_{k}^1 \\ 0 \end{array}\right),
\mathbf{K}_{k}=\left(\begin{array}{cc} 0 \\ (H_{k}^1)^TR_2H^2_{k} \\[1mm] C_{k}^T \end{array}\right),
\mathbf{G}=\left(\begin{array}{ccc} G_2 \\ 0 \\ 0 \end{array}\right), \\
\mathbf{L}_{k}=\left(\begin{array}{ccc} A^T & 0 & 0 \\ -(H_{k}^1)^TB_1^T & \widetilde{A}_{k}^T & (H_{k}^1)^TR_2H^3_{k} \\ [1mm] 0 & 0 & \widetilde{A}_{k}^T \end{array}\right),
\widetilde{\mathbf{C}}_k=\left(\begin{array}{ccc} 0 \\ 0 \\ \widetilde{C}_k^T \end{array}\right),
\mathbf{S}_k=\left(\begin{array}{ccc} 0 \\ 0 \\ (H^3_k)^T \end{array}\right), \\
\mathbf{D}_k=\left(\begin{array}{ccc} \widetilde{B}_k-B_1H_k^1\sum_{i=t}^{k-1}(D_i^{(k)})^T, \\[2mm]
\widetilde{B}_k+\widetilde{A}_k\sum_{i=t}^{k-1}(D_i^{(k)})^T \\[2mm] (H^3_k)^TR_2H^2_k+(H^3_k)^TR_2H_k^1\sum_{i=t}^{k-1}(D_i^{(k)})^T \end{array}\right),
\end{array}
\right.
\end{eqnarray}
then the FB$\Delta$Es  in (\ref{backwards equation-21}), (\ref{state-X*}), (\ref{backwards equation-2-211}) are equivalently rewritten as
\begin{eqnarray}
\label{FBDE-2}
\left\{
\begin{array}{l}
X_{k+1}^*=\widetilde{A}_{k}X_{k}^*+\widetilde{B}_{k}v_{k}^*+\widetilde{\mathbf{C}}_{k}^T\mathbf{Z}_{k+1}^{*}, \\[1mm]
\mathbf{Z}_{k}^{*}=\mathbf{H}_{k}X_{k}^*+\mathbf{K}_{k}v_{k}^*
+\mathbf{L}_{k}\mathbf{Z}_{k+1}^{*}, \\[1mm]
X_t^*=x, ~~~\mathbf{Z}_{N}^{*}=\mathbf{G}X_N^*, ~~~k\in \mathbb{T}_t,
\end{array}
\right.
\end{eqnarray}
and the stationary condition (\ref{stationary-condition-1}) is equivalent to
\begin{eqnarray}
\label{stationary-condition-2}
&&0=\Big{[}W_2+(H^2_k)^TR_2H^2_k+\sum_{i=t}^{k-1}D_i^{(k)}(H_k^1)^TR_2H^2_k\Big{]}v_k^* \nonumber\\
&&\hphantom{0=}+\Big{[}(H^2_k)^TR_2H_k^1
+\sum_{i=t}^{k-1}D_i^{(k)}(H_k^1)^TR_2H^1_k\Big{]}X_k^*+\mathbf{D}_k^T\mathbf{Z}_{k+1}^{*},  ~~~k\in \mathbb{T}_t.
\end{eqnarray}
Letting
\begin{eqnarray}
\left\{
\begin{array}{l}
F_k=W_2+(H^2_k)^TR_2H^2_k+\sum_{i=t}^{k-1}D_i^{(k)}(H_k^1)^TR_2H^2_k, \\[1mm]
O_k=(H^2_k)^TR_2H_k^1+\sum_{i=t}^{k-1}D_i^{(k)}(H_k^1)^TR_2H^1_k, \\[1mm]
k\in\mathbb{T}_t,
\end{array}
\right.
\end{eqnarray}
then (\ref{stationary-condition-2}) is denoted as
\begin{eqnarray}
\label{stationary-condition-3}
0=F_kv_k^*+O_kX_k^*+\mathbf{D}_k^T\mathbf{Z}_{k+1}^{*},  ~~~k\in \mathbb{T}_t.
\end{eqnarray}

{The following theorem characterizes the unique existence of $v^*\in l^2(\mathbb{T}_t; \mathbb{R}^{m_2})$ in (\ref{inequality-2}), by virtue of some matrices with the nonsingularity constraint, which can be easily checked.
}

\begin{theorem}\label{Theorem-stationary-condition}

The following statements are equivalent.
\begin{itemize}


\item[i)]  There exists a unique $v^*\in l^2(\mathbb{T}_t; \mathbb{R}^{m_2})$ such that (\ref{inequality-2}) holds for any $k\in\mathbb{T}_t$ and any $v_k\in l^2(k; \mathbb{R}^{m_2})$.

%

\item[ii)] $F_k$ and $\mathbf{I}-(\widetilde{\mathbf{C}}_{k}^T-\widetilde{B}_{k}F_{k}^{-1}\mathbf{D}_{k}^T)\mathbf{T}_{k+1}, k\in \mathbb{T}_t$ are invertible, where
\begin{eqnarray}
\label{T}
\left\{
\begin{array}{l}
\mathbf{T}_k=(\mathbf{L}_{k}-\mathbf{K}_{k}F_{k}^{-1}\mathbf{D}_k^T)\mathbf{T}_{k+1}\big{[}\mathbf{I}-(\widetilde{\mathbf{C}}_{k}^T-\widetilde{B}_{k}F_{k}^{-1}\mathbf{D}_k^T)
\mathbf{T}_{k+1}\big{]}^{-1}(\widetilde{A}_{k}-\widetilde{B}_{k}F_{k}^{-1}O_{k}) \\[1mm]
\hphantom{\mathbf{T}_k=}+\mathbf{H}_{k}-\mathbf{K}_{k}F_{k}^{-1}O_{k}, \\[1mm]
\mathbf{T}_N=\mathbf{G}, ~~~k\in \mathbb{T}_t.
\end{array}
\right.
\end{eqnarray}

\end{itemize}

In this case, the $FB\Delta E$ of (\ref{FBDE-2}) has the following expression
\begin{eqnarray}
\label{FBDE-3}
\left\{
\begin{array}{l}
X_{k+1}^*=\big{[}\mathbf{I}-(\widetilde{\mathbf{C}}_{k}^T-\widetilde{B}_{k}F_{k}^{-1}\mathbf{D}_k^T)\mathbf{T}_{k+1}\big{]}^{-1}(\widetilde{A}_{k}
-\widetilde{B}_{k}F_{k}^{-1}O_{k})X_{k}^*, \\[1mm]
\mathbf{Z}_{k}^{*}=\mathbf{T}_{k}X_{k}^*, \\[1mm]
X_t^*=x, ~~~\mathbf{Z}_{N}^{*}=\mathbf{G}X_N^*, ~~~k\in \mathbb{T}_t,
\end{array}
\right.
\end{eqnarray}
and $v^*_k$ of (\ref{stationary-condition-3}) and i) is computed via
\begin{eqnarray}
\label{leader's optimal control-3}
\hspace{-1em}v_k^*=-F_k^{-1}\Big{\{}O_k+\mathbf{D}_k^T\mathbf{T}_{k+1}\big{[}\mathbf{I}-(\widetilde{\mathbf{C}}_{k}^T
-\widetilde{B}_{k}F_{k}^{-1}\mathbf{D}_k^T)\mathbf{T}_{k+1}\big{]}^{-1}(\widetilde{A}_{k}
-\widetilde{B}_{k}F_{k}^{-1}O_{k})\Big{\}}X_k^*, ~~k\in \mathbb{T}_t.
\end{eqnarray}

\end{theorem}

We finally characterize the control process $u^*$ of $iii)$ of Definition \ref{definition}. Substituting the $v$ of (\ref{follower's optimal control}) with (\ref{leader's optimal control-3}) and noting (\ref{FBDE-3}), the unique $u^*\in l^2(\mathbb{T}_t; \mathbb{R}^{m_1})$ is expressed as
\begin{eqnarray}\label{u*}
&&\hspace{-2em}u_k^*=[\alpha^t(x, v^*)]_k=-\big{(}H^1_kX_k^*+H^2_kv_k^*+H^3_k\pi_{k+1}^*\big{)} \nonumber\\
%
%
%
&&\hspace{-2em}\hphantom{u^*}=\Big{\{}(H^2_kF_k^{-1}\mathbf{D}_k^T-\mathbf{S}_k^T)\mathbf{T}_{k+1}\big{[}\mathbf{I}-(\widetilde{\mathbf{C}}_{k}^T
-\widetilde{B}_{k}F_{k}^{-1}\mathbf{D}_k^T)\mathbf{T}_{k+1}\big{]}^{-1}(\widetilde{A}_{k}
-\widetilde{B}_{k}F_{k}^{-1}O_{k}) \nonumber\\
&&\hspace{-2em}\hphantom{u^*=\Big{\{}}-H^1_k+H^2_kF_k^{-1}O_k\Big{\}}X_{k}^*,
 ~~~~~~k\in \mathbb{T}_t.
\end{eqnarray}

To conclude this section, we have the following result.

\begin{theorem}\label{Theorem}

The following statements are equivalent.
\begin{itemize}

\item[i)] The Stackelberg game associated with (\ref{cost-1}) (\ref{system-dynamics}) admits a unique open-loop  equilibrium solution $(u^*, v^*)$ for the initial pair $(t,x)$.

\item[ii)] $M_k, k\in \mathbb{T}_t$, of (\ref{M}) are positive definite, and $F_k$,  $\mathbf{I}-(\widetilde{\mathbf{C}}_{k}^T-\widetilde{B}_{k}F_{k}^{-1}\mathbf{D}_{k}^T)\mathbf{T}_{k+1}, k\in \mathbb{T}_t$ are invertible.

\end{itemize}

In this case, the controls $u^*, v^*$ of i) are given by (\ref{u*}) and (\ref{leader's optimal control-3}), respectively.

\end{theorem}

\subsection{About time consistetncy}

This section is about the time consistency of  open-loop equilibrium solution of Definition \ref{definition}. Let $(u^*, v^*)$ be the unique open-loop equilibrium solution for the initial pair $(t,x)$. Along the equilibrium system (\ref{system-dynamics-3}), at time instant $\tau$ we reconsider the Stackelberg game, and denote its unique open-loop  equilibrium solution (for the initial pair $(\tau, X^*_\tau)$), if it exists,  as $(u^{\tau*}, v^{\tau *})$, which is similarly defined as that of Definition \ref{definition}:
\begin{itemize}
\item[i)] There exists a unique map $\alpha^\tau: \mathbb{R}^n\times l^2(\mathbb{T}_\tau; \mathbb{R}^{m_2})\mapsto l^2(\mathbb{T}_\tau; \mathbb{R}^{m_1})$ such that
    \begin{eqnarray*}\label{inequality-1-c}
      J_1(k, X_k; {\alpha^\tau}(x, v^\tau)|_{\mathbb{T}_k}, v^\tau|_{\mathbb{T}_k})\leq J_1(k, X_k; (u_k, {\alpha^\tau}(x, v^\tau)|_{\mathbb{T}_{k+1}}), v^\tau|_{\mathbb{T}_k})
    \end{eqnarray*}
    holds for each $v^\tau\in l^2(\mathbb{T}_\tau; \mathbb{R}^{m_2})$, any $k\in\mathbb{T}_\tau$ and any $u_k\in l^2(k; \mathbb{R}^{m_1})$, where the state $X_k$ is computed via
    \begin{eqnarray*}\label{system-dynamics-2-c}
\left\{
\begin{array}{l}
X_{k+1}=AX_k+B_1[\alpha^\tau(x, v^\tau)]_{k}+B_2v^\tau_k,\\
X_\tau=X_\tau^*,~~~~k\in \mathbb{T}_\tau.
\end{array}
\right.
\end{eqnarray*}

\item[ii)]  There exists a unique $v^{\tau*}\in l^2(\mathbb{T}_\tau; \mathbb{R}^{m_2})$ such that
    \begin{eqnarray}\label{inequality-2-c}
    J_2(k, X_k^{\tau*}; {\alpha^\tau}(x, v^{\tau*})|_{\mathbb{T}_k}, v^{\tau*}|_{\mathbb{T}_k})\leq J_2(k, X_k^{\tau*}; {\alpha^\tau}(x, v^{\tau*-k})|_{\mathbb{T}_k}, (v_k, v^{\tau*}|_{\mathbb{T}_{k+1}}))
    \end{eqnarray}
holds for any $k\in\mathbb{T}_\tau$ and any $v_k\in l^2(k; \mathbb{R}^{m_2})$, %
where
\begin{eqnarray*}
v^{\tau*-k}_\ell=\left\{
\begin{array}{ll}
v_k,  ~~&\ell=k,\\
v^*_\ell,  ~~&\ell\neq k, \ell \in \mathbb{T}_\tau.
\end{array}
\right.
\end{eqnarray*}

The state $X_k^{\tau*}$ in (\ref{inequality-2-c}) is computed via
\begin{eqnarray}\label{system-dynamics-3-c}
\left\{
\begin{array}{l}
X^{\tau*}_{k+1}=AX^{\tau*}_k+B_1[\alpha^\tau(x, v^{\tau*})]_{k}+B_2v^{\tau*}_k,\\
X^{\tau*}_\tau=X^*_\tau,~~~~k\in \mathbb{T}_\tau,
\end{array}
\right.
\end{eqnarray}

\item[iii)]  $u^{\tau *}=\alpha^\tau(x, v^{\tau*})$.

\end{itemize}

By mimicing all the derivations of above section, we have a result on open-loop equilibrium solution for the initial pair $(\tau, X^*_\tau)$ that is parallel to Theorem \ref{Theorem}.

\begin{theorem}\label{Theorem-c}

The following statements are equivalent.
\begin{itemize}

\item[i)] The Stackelberg game associated with (\ref{cost-1}) (\ref{system-dynamics}) admits a unique open-loop  equilibrium solution $(u^{\tau*}, v^{\tau*})$ for the initial pair $(\tau, X^*_\tau)$.

\item[ii)] $M_k, k\in \mathbb{T}_\tau$, of (\ref{M}) are positive definite, and $F_k$, $\mathbf{I}-(\widetilde{\mathbf{C}}_{k}^T-\widetilde{B}_{k}F_{k}^{-1}\mathbf{D}_{k}^T)\mathbf{T}_{k+1}, k\in \mathbb{T}_\tau$ are invertible.

\end{itemize}

In this case, the controls $u^{\tau *}, v^{\tau *}$ of i) are given by
\begin{eqnarray*}\label{u*c}
&&\hspace{-3em}u_k^{\tau *}=\Big{\{}(H^2_kF_k^{-1}\mathbf{D}_k^T-\mathbf{S}_k^T)\mathbf{T}_{k+1}\big{[}\mathbf{I}-(\widetilde{\mathbf{C}}_{k}^T
-\widetilde{B}_{k}F_{k}^{-1}\mathbf{D}_k^T)\mathbf{T}_{k+1}\big{]}^{-1}(\widetilde{A}_{k}
-\widetilde{B}_{k}F_{k}^{-1}O_{k}) \nonumber\\
&&\hspace{-3em}\hphantom{u^{\tau *}=\Big{\{}}-H^1_k+H^2_kF_k^{-1}O_k\Big{\}}X_{k}^{\tau*},\\
&&\hspace{-3em}\label{v^*c}v_k^{\tau *}=-F_k^{-1}\Big{\{}O_k+\mathbf{D}_k^T\mathbf{T}_{k+1}\big{[}\mathbf{I}-(\widetilde{\mathbf{C}}_{k}^T
-\widetilde{B}_{k}F_{k}^{-1}\mathbf{D}_k^T)\mathbf{T}_{k+1}\big{]}^{-1}(\widetilde{A}_{k}
-\widetilde{B}_{k}F_{k}^{-1}O_{k})\Big{\}}X_k^{\tau*}, ~~k\in \mathbb{T}_\tau
\end{eqnarray*}
with
\begin{eqnarray*}
\label{FBDE-3-2}
\left\{
\begin{array}{l}
X_{k+1}^{\tau*}=\big{[}\mathbf{I}-(\widetilde{\mathbf{C}}_{k}^T-\widetilde{B}_{k}F_{k}^{-1}\mathbf{D}_k^T)\mathbf{T}_{k+1}\big{]}^{-1}(\widetilde{A}_{k}
-\widetilde{B}_{k}F_{k}^{-1}O_{k})X_{k}^{\tau*}, \\[1mm]
X_\tau^{\tau*}=X_\tau^*, ~~~k\in \mathbb{T}_\tau,
\end{array}
\right.
\end{eqnarray*}

\end{theorem}

Comparing Theorem \ref{Theorem} and Theorem \ref{Theorem-c} and noticing the expressions of $X^*, X^{\tau*}$, $(u^*, v^*), (u^{\tau*}$, $v^{\tau*})$, we must have the following facts.

\begin{theorem}
The following facts hold.

\begin{itemize}

\item[i)] If the Stackelberg game associated with (\ref{cost-1}) (\ref{system-dynamics}) admits a unique open-loop  equilibrium solution $(u^*, v^*)$ for the initial pair $(t, x)$, then for any equilibrium pair $(\tau, X^{*}_\tau)$ with $\tau\in \mathbb{T}_t$ the Stackelberg game admits a unique open-loop equilibrium solution $(u^{\tau*}, v^{\tau*})$.

\item[ii)] The open-loop equilibrium solution is time-consistent, namely, for $\tau\in \mathbb{T}_t$, $u^*|_{\mathbb{T}_\tau}=u^{\tau*}, v^*|_{\mathbb{T}_\tau}=v^{\tau*}$ with the control inputs given in i).

\end{itemize}

\end{theorem}

\begin{remark}
The initial state of  (\ref{system-dynamics-3-c}) is $X^*_\tau$, which means that on $\{t,...,\tau-1\}$ we adopt the equilibrium control $(\{u^*_t,..., u^*_{\tau-1}\}, \{v^*_t,..., v^*_{\tau-1}\})$. Hence, the time consistency of open-loop equilibrium solution is weak in the sense of \cite{Basar}.
We now show the reason that ensure the time consistency of open-loop equilibrium solution. Note that under the follower's best response the controlled system of the leader is (\ref{system-leader}), which is decoupled, namely, forward state does not appear in the backward difference equation. Then, the perturbation $v^\varepsilon$ of Proposition \ref{variation} cannot influence the backward difference equation, though it will influence the forward difference equation. Hence, we need not introduce a adjoint forward difference equation for the backward difference equation; this is why in (\ref{FBDE-2}) the forward state is only $X^*$, i.e., we need not enlarge the forward state. On the other hand, in characterizing the open-loop Stackelberg solution, we need to enlarge the equilibrium state $\widehat{X}_k$ to $\xi^{(0,x)}_k=((\beta_k^{(0,x)})^T, \widehat{X}_k^T)^T$ of (\ref{closed-loop-system-0}). As shown in Section \ref{Section-2-1}, it is the term $\beta^{(0,x)}$ that ruins the time consistency of open-loop Stackelberg solution.

\end{remark}

\section{Proofs}\label{proof}

\subsection{Proof of Proposition \ref{variation}}

Under $v^\varepsilon$, the follower's control input is
\begin{eqnarray}
[\alpha^t(x, v^\varepsilon)]_k=-(H_k^1X^\varepsilon_k+H^2_kv^{\varepsilon}_k+H_k^3\pi_{k+1}^\varepsilon),~~~k\in \mathbb{T}_t
\end{eqnarray}
with $X^\varepsilon$ given in (\ref{system-leader-e}).
%
%
%
Hence, one gets
\begin{eqnarray*}
\left\{\begin{array}{l}
\frac{X_{\ell+1}^{\varepsilon}-X_{\ell+1}}{\varepsilon}=\widetilde{A}_\ell\frac{X_{\ell}^{\varepsilon}
-X_{\ell}}{\varepsilon},\\[1mm]
\frac{X_{k+1}^{\varepsilon}-X_{k+1}}{\varepsilon}=\widetilde{A}_k\frac{X_{k}^{\varepsilon}
-X_{k}}{\varepsilon}+\widetilde{B}_k\widetilde{v}_k,\\[1mm]
\frac{X_{k}^{\varepsilon}-X_{k}}{\varepsilon}=\widetilde{A}_{k-1}\frac{X_{k-1}^\varepsilon
-X_{k-1}}{\varepsilon}+\widetilde{C}_{k-1}C_k^T\widetilde{v}_k,\\[1mm]
\frac{X_{i+1}^{\varepsilon}-X_{i+1}}{\varepsilon}=\widetilde{A}_{i}\frac{X_{i}^{\varepsilon}
-X_{i}}{\varepsilon}+\widetilde{C}_{i}
\widetilde{A}_{i+1}^T\widetilde{A}_{i+2}^T\ldots\widetilde{A}_{k-1}^TC_{k}^T\widetilde{v}_k,\\[1mm]
\frac{X_{t}^{\varepsilon}-X_{t}}{\varepsilon}=0,\\[1mm]
\ell\in\mathbb{T}_{k+1},~~i\in {\{}t, t+1, \ldots, k-2{\}}.
\end{array}
\right.
\end{eqnarray*}
Denoting $\frac{X_\ell^{\varepsilon}-X_\ell}{\varepsilon}$ by $\eta^{(k)}_\ell$, we have (\ref{Y-1}) and $X_\ell^{\varepsilon}=X_\ell+\varepsilon \eta^{(k)}_\ell$, $\ell \in \mathbb{T}_t$.
It should be noted that the internal state of $J_2(k, X_k; \alpha^t(x, v^\varepsilon)|_{\mathbb{T}_k}, (v_k+\varepsilon\widetilde{v}_k,  v|_{\mathbb{T}_{k+1}}))$ is
\begin{eqnarray}
\left\{
\begin{array}{l}
X_{\ell+1}^c=AX_\ell^c+B_1[\alpha^t(x, v^{\varepsilon})]_{\ell}+B_2v^{\varepsilon}_\ell,\\
X_k^c=X_k,~~\ell \in \mathbb{T}_k,
\end{array}
\right.
\end{eqnarray}
which is different from $\{X_{\ell}^\varepsilon, \ell\in \mathbb{T}_k\}$.
Furthermore, one has
\begin{eqnarray*}
\left\{
\begin{array}{l}
\frac{X_{\ell+1}^c-X_{\ell+1}}{\varepsilon}=A\frac{X_{\ell}^c-X_{\ell}}{\varepsilon}
-B_1H_{\ell}^1\frac{X_{\ell}^\varepsilon-X_{\ell}}{\varepsilon},\\[1mm]
\frac{X_{k+1}^c-X_{k+1}}{\varepsilon}=A\frac{X_{k}^c-X_{k}}{\varepsilon}
-B_1H_{k}^1\frac{X_{k}^\varepsilon-X_{k}}{\varepsilon}+\widetilde{B}_k\widetilde{v}_k,\\[1mm]
\frac{X_{k}^c-X_{k}}{\varepsilon}=0, ~~~\ell\in\mathbb{T}_{k+1}.
\end{array}
\right.
\end{eqnarray*}
Denoting $\frac{X_\ell^{c}-X_\ell}{\varepsilon}$ by $\xi_\ell$, we have (\ref{Y-2}) and $X_\ell^{c}=X_\ell+\varepsilon\xi_\ell$, $\ell\in\mathbb{T}_k$. Noting that $\pi_\ell$, $\ell\in\mathbb{T}_{k+1}$ are not influenced by $v_k+\varepsilon\widetilde{v}_k$, i.e., $\pi_\ell^\varepsilon=\pi_\ell, \ell\in \mathbb{T}_{k+1}$. Then, we have
\begin{eqnarray}
\label{difference formula-3}
&&\hspace{-1.85em}J_2(k, X_k; \alpha^t(x, v^\varepsilon)|_{\mathbb{T}_{k}}, (v_k+\varepsilon \widetilde{v}_k,  v|_{\mathbb{T}_{k+1}}))-
J_2(k, X_k; \alpha^t(x, v)|_{\mathbb{T}_{k}}, v|_{\mathbb{T}_{k}})\nonumber\\
&&\hspace{-1.85em}=\varepsilon^2\Big{[}\sum_{\ell=k}^{N-1}\xi_\ell^TQ_2\xi_\ell+\widetilde{v}_k^TW_{2}\widetilde{v}_k
+\xi_{N}^TG_2\xi_{N}+(H_k^1\eta^{(k)}_k+H^2_k\widetilde{v}_k)^TR_2(H_k^1\eta^{(k)}_k+H^2_k\widetilde{v}_k)\nonumber\\
&&\hspace{-1.85em}\hphantom{=}+
\sum_{\ell=k+1}^{N-1}(\eta_{\ell}^{(k)})^T(H_{\ell}^1)^TR_2H_{\ell}^1\eta^{(k)}_{\ell}\Big{]}+2\varepsilon\Big{\{}\sum_{\ell=k}^{N-1}\Big{[}X_\ell^TQ_2\xi_{\ell}+(H_{\ell}^1X_{\ell}+H^2_{\ell} v_{\ell}+H_{\ell}^3\pi_{\ell+1})^TR_2H_{\ell}^1\eta^{(k)}_\ell\Big{]}\nonumber\\
&&\hspace{-1.85em}\hphantom{=}+v_k^TW_2\widetilde{v}_k+X_{N}^TG_2\xi_{N}+(H_k^1X_k+H^2_kv_k+H_k^3\pi_{k+1})^T R_2H^2_k\widetilde{v}_k\Big{\}}.
\end{eqnarray}
As $\xi_k=0$, $\eta^{(k)}_t=0$, it holds that
\begin{eqnarray*}
&&\sum_{\ell=k}^{N-1}\Big{[}X_\ell^TQ_2\xi_{\ell}+(H_{\ell}^1X_{\ell}+H^2_{\ell} v_{\ell}+H_{\ell}^3\pi_{\ell+1})^TR_2H_{\ell}^1\eta^{(k)}_\ell\Big{]}+v_k^TW_2\widetilde{v}_k\\
&&+X_{N}^TG_2\xi_{N}+(H_k^1X_k+H^2_kv_k+H_k^3\pi_{k+1})^TR_2H^2_k\widetilde{v}_k\\
&&=\sum_{\ell=k}^{N-1}\Big{[}X_{\ell}^TQ_2\xi_{\ell}+(H_{\ell}^1X_{\ell}+H^2_{\ell}v_{\ell}
+H_{\ell}^3\pi_{\ell+1})^TR_2H_{\ell}^1\eta^{(k)}_{\ell}+Z_{\ell+1}^T\xi_{\ell+1}-Z_\ell^T\xi_{\ell}\\
&&\hphantom{=}
+(\overline{Z}_{\ell+1}^{(k)})^T\eta^{(k)}_{\ell+1}-(\overline{Z}_{\ell}^{(k)})^T\eta^{(k)}_{\ell}\Big{]}+\sum_{i=t}^{k-1}
\Big{[}(\overline{Z}_{i+1}^{(k)})^T\eta^{(k)}_{i+1}-(\overline{Z}_{i}^{(k)})^T\eta^{(k)}_{i}\Big{]}+v_k^TW_2\widetilde{v}_k\\
&&\hphantom{=}+(H_k^1X_k+H^2_kv_k+H_k^3\pi_{k+1})^TR_2H^2_k\widetilde{v}_k\\
%
%
&&=\sum_{\ell=k}^{N-1}\Big{\{}(Q_2X_{\ell}+A^TZ_{\ell+1}-Z_{\ell})^T\xi_{\ell}+\Big{[}(H_{\ell}^1)^TR_2
(H_{\ell}^1X_{\ell}+H^2_{\ell}v_{\ell}+H_{\ell}^3\pi_{\ell+1})\\
&&\hphantom{=}-(H_{\ell}^1)^TB_1^TZ_{\ell+1}+\widetilde{A}_{\ell}^T\overline{Z}_{\ell+1}^{(k)}-
\overline{Z}_{\ell}^{(k)}\Big{]}^T\eta^{(k)}_{\ell}\Big{\}}+\sum_{i=t}^{k-1}
(\widetilde{A}_{i}^T\overline{Z}_{i+1}^{(k)}-\overline{Z}_{i}^{(k)})^T\eta^{(k)}_{i}\\
&&\hphantom{=}+\Big{[}\widetilde{B}_k^TZ_{k+1}+\widetilde{B}_{k}^T\overline{Z}_{k+1}^{(k)}+W_2v_k
+(H^2_k)^TR_2(H_k^1X_k+H^2_kv_k+H_k^3\pi_{k+1})\\
&&\hphantom{=}+\sum_{i=t}^{k-1}C_{k}\widetilde{A}_{k-1}\cdots \widetilde{A}_{i+1}\widetilde{C}_{i}^T\overline{Z}_{i+1}^{(k)}
\Big{]}^T\widetilde{v}_k.
\end{eqnarray*}
This and (\ref{difference formula-3}) imply the conclusion, and the proof is completed.

\subsection{Proof of Theorem \ref{equivalence-1}}

$i)\Rightarrow ii)$. Combing Proposition \ref{variation} and Lemma \ref{stationary condition}, one gets
\begin{eqnarray}
\label{difference-formula-4}
&&J_2(k, X_k^*; \alpha^t(x, v^{*-k})|_{\mathbb{T}_{k}}, (v_k^*+\varepsilon \widetilde{v}_k,  v^*|_{\mathbb{T}_{k+1}}))-J_2(k, X_k^*; \alpha^t(x, v^*)|_{\mathbb{T}_{k}}, v^*|_{\mathbb{T}_{k}})\nonumber\\
&&=2\varepsilon \Big{\{}\Big{[}\widetilde{B}_{k}^T-\sum_{i=t}^{k-1}D_i^{(k)}(H_k^1)^TB_1^T\Big{]}Z_{k+1}^*+(\widetilde{B}_{k}^T
+\sum_{i=t}^{k-1}D_i^{(k)}\widetilde{A}_k^T)\overline{Z}_{k+1}^{(k)*} \nonumber\\
&&\hphantom{=2\varepsilon\Big{[}}+\Big{[}(H^2_k)^TR_2H^1_k+\sum_{i=t}^{k-1}D_i^{(k)}(H_k^1)^TR_2H_k^1\Big{]}X_k^* \nonumber\\
&&\hphantom{=2\varepsilon\Big{[}}+\Big{[}W_2+(H^2_k)^TR_2H^2_k+\sum_{i=t}^{k-1}D_i^{(k)}(H_k^1)^TR_2H^2_k\Big{]}v_k^*
 \nonumber\\
&&\hphantom{=2\varepsilon\Big{[}}+\Big{[}(H^2_k)^TR_2H_k^3+\sum_{i=t}^{k-1}D_i^{(k)}(H_k^1)^TR_2H_k^3\Big{]}\pi_{k+1}^* \Big{\}}^T\widetilde{v}_k +\varepsilon^2\widehat{J}_2(k, 0; \widetilde{v}_k) \nonumber\\
&&\geq 0.
\end{eqnarray}
It can be seen from (\ref{convex_2}) that $\widehat{J}_2(k, 0; \widetilde{v}_k)\geq0$ always holds for any $\widetilde{v}_k \in l^2(k; \mathbb{R}^{m_2})$. As (\ref{difference-formula-4}) holds for any $\varepsilon \in \mathbb{R}$ and any $\widetilde{v}_k \in l^2(k; \mathbb{R}^{m_2})$, one must obtain (\ref{stationary-condition-1}). In fact,
if for some $k_1 \in \mathbb{T}_t$,
\begin{eqnarray*}
&&\sigma_{k_1}\triangleq\Big{[}\widetilde{B}_{k_1}^T-\sum_{i=t}^{k_1-1}D_i^{(k_1)}(H_{k_1}^1)^TB_1^T\Big{]}Z_{k_1+1}^*+(\widetilde{B}_{k_1}^T
+\sum_{i=t}^{k_1-1}D_i^{(k_1)}\widetilde{A}_{k_1}^T)\overline{Z}_{k_1+1}^{(k_1)*} \\
&&\hphantom{\sigma_{k_1}=}+\Big{[}(H^2_{k_1})^TR_2H_{k_1}^1
+\sum_{i=t}^{k_1-1}D_i^{(k_1)}(H_{k_1}^1)^TR_2H_{k_1}^1\Big{]}X_{k_1}^*  \\
&&\hphantom{\sigma_{k_1}=}+\Big{[}W_2+(H^2_{k_1})^TR_2H^2_{k_1}+\sum_{i=t}^{k_1-1}D_i^{(k_1)}(H_{k_1}^1)^TR_2H^2_{k_1}\Big{]}v_{k_1}^*\\
&&\hphantom{\sigma_{k_1}=}+\Big{[}(H^2_{k_1})^TR_2H_{k_1}^3+\sum_{i=t}^{k_1-1}D_i^{(k_1)}(H_{k_1}^1)^TR_2H_{k_1}^3\Big{]}\pi_{k_1+1}^* \\
&&\hphantom{\sigma_{k_1}}\neq 0,
\end{eqnarray*}
let $\widetilde{v}_{k_1}=\sigma_{k_1}$; this and (\ref{difference-formula-4}) imply that
\begin{eqnarray}
\label{Thm-leader-proof-1}
\varepsilon^2\widehat{J}_2(k_1, 0; \sigma_{k_1})+2\varepsilon|\sigma_{k_1}|^2\geq 0
\end{eqnarray}
holds for any $\varepsilon \in \mathbb{R}$. 
If $\varepsilon$ is a negative number with $|\varepsilon|$ sufficiently small, one has
\begin{eqnarray*}
\varepsilon^2\widehat{J}_2(k, 0; \sigma_{k_1})+2\varepsilon|\sigma_{k_1}|^2< 0,
\end{eqnarray*}
which contradicts (\ref{Thm-leader-proof-1}). Therefore, $\sigma_{k_1}$ must be 0, and (\ref{stationary-condition-1}) holds.

$ii)\Rightarrow i)$. In this case, one get
\begin{eqnarray*}
J_2(k, X_k^*; {\alpha^t}(x, v^{*-k})|_{\mathbb{T}_k}, (v_k, v^*|_{\mathbb{T}_{k+1}}))-
J_2(k, X_k^*; {\alpha^t}(x, v^*)|_{\mathbb{T}_k}, v^*|_{\mathbb{T}_k})\geq 0,
\end{eqnarray*}
and the conclusion follows.

\subsection{Proof of Theorem \ref{Theorem-stationary-condition}}

\emph{i)$\Rightarrow$ii)}. By Theorem \ref{equivalence-1}, let $v^*\in l^2(\mathbb{T}_t; \mathbb{R}^{m_2})$ be the one such that the stationary condition (\ref{stationary-condition-1}) holds. Under this $v^*$, equations (\ref{backwards equation-21}), (\ref{state-X*}) and (\ref{backwards equation-2-211}) have the unique solution $Z_k^*$, $(\pi_k^*, X_k^*)$ and $\overline{Z}_k^*$, respectively. Equivalently, (\ref{stationary-condition-3}) has a unique solution $(X^*, \mathbf{Z}^*)$ and (\ref{stationary-condition-3}) holds. It is easy to get that $F_k$, $k\in\mathbb{T}_t$ in (\ref{stationary-condition-3}) is invertible due to the uniqueness of $v^*$.
Hence, (\ref{stationary-condition-3}) is equal to
\begin{eqnarray*}
v_k^*=-F_k^{-1}(O_kX_k^*+\mathbf{D}_k^T\mathbf{Z}_{k+1}^{*}), ~~~k\in \mathbb{T}_t.
\end{eqnarray*}
Accordingly, the FB$\Delta$E (\ref{FBDE-2}) is rewritten as
\begin{eqnarray}
\label{FBDE-5}
\left\{
\begin{array}{l}
X_{k+1}^*=(\widetilde{A}_{k}-\widetilde{B}_{k}F_k^{-1}O_k)X_{k}^*+(\widetilde{\mathbf{C}}_{k}^T-\widetilde{B}_{k}F_k^{-1}\mathbf{D}_k^T)\mathbf{Z}_{k+1}^{*}, \\[1mm]
\mathbf{Z}_{k}^{*}=(\mathbf{H}_{k}-\mathbf{K}_{k}F_k^{-1}O_k)X_{k}^*
+(\mathbf{L}_{k}-\mathbf{K}_{k}F_k^{-1}\mathbf{D}_k^T)\mathbf{Z}_{k+1}^{*}, \\[1mm]
X_t^*=x, ~~~\mathbf{Z}_{N}^{*}=\mathbf{G}X_N^*, ~~~k\in \mathbb{T}_t,
\end{array}
\right.
\end{eqnarray}
which is uniquely solvable.

Noting the terminal condition $\mathbf{Z}_{N}^{*}=\mathbf{G}X_N^*$ and the first equation of (\ref{FBDE-5}), one has
\begin{eqnarray*}
X_N^*=(\widetilde{A}_{N-1}-\widetilde{B}_{N-1}F_{N-1}^{-1}O_{N-1})X_{N-1}^*
+(\widetilde{\mathbf{C}}_{N-1}^T-\widetilde{B}_{N-1}F_{N-1}^{-1}\mathbf{D}_{N-1}^T)\mathbf{G}X_N^*,
\end{eqnarray*}
that is,
\begin{eqnarray*}
\big{[}\mathbf{I}-(\widetilde{\mathbf{C}}_{N-1}^T-\widetilde{B}_{N-1}F_{N-1}^{-1}\mathbf{D}_{N-1}^T)\mathbf{G}\big{]}X_N^*
=(\widetilde{A}_{N-1}-\widetilde{B}_{N-1}F_{N-1}^{-1}O_{N-1})X_{N-1}^*,
\end{eqnarray*}
Using the unique solvability of (\ref{FBDE-5}), it yields that $\mathbf{I}-(\widetilde{\mathbf{C}}_{N-1}^T-\widetilde{B}_{N-1}F_{N-1}^{-1}\mathbf{D}_{N-1}^T)\mathbf{G}$ is invertible. It then follows that $\mathbf{Z}_{N-1}^{*}=\mathbf{T}_{N-1}X_{N-1}^*$ with $\mathbf{T}_{N-1}$ satisfying (\ref{T}) for  $k=N-1$.

Assume now that $\mathbf{Z}_{k+1}^{*}=\mathbf{T}_{k+1}X_{k+1}^*$ holds. Let us show $\mathbf{Z}_{k}^{*}=\mathbf{T}_{k}X_{k}^*$. By substituting $\mathbf{Z}_{k+1}^{*}=\mathbf{T}_{k+1}X_{k+1}^*$ into (\ref{FBDE-5}), we have
\begin{eqnarray*}
X_{k+1}^*=(\widetilde{A}_{k}-\widetilde{B}_{k}F_k^{-1}O_k)X_{k}^*+(\widetilde{\mathbf{C}}_{k}^T-\widetilde{B}_{k}F_k^{-1}\mathbf{D}_k^T)\mathbf{T}_{k+1}X_{k+1}^*,
\end{eqnarray*}
namely,
\begin{eqnarray*}
\big{[}\mathbf{I}-(\widetilde{\mathbf{C}}_{k}^T-\widetilde{B}_{k}F_{k}^{-1}\mathbf{D}_{k}^T)\mathbf{T}_{k+1}\big{]}X_{k+1}^*
=(\widetilde{A}_{k}-\widetilde{B}_{k}F_{k}^{-1}O_{k})X_{k}^*,
\end{eqnarray*}
Using the unique solvability of (\ref{FBDE-5}), it yields that $\mathbf{I}-(\widetilde{\mathbf{C}}_{k}^T-\widetilde{B}_{k}F_{k}^{-1}\mathbf{D}_{k}^T)\mathbf{T}_{k+1}$ is invertible. Combining this with (\ref{FBDE-5}), it yields that
\begin{eqnarray*}
\mathbf{Z}_{k}^{*}=\mathbf{T}_{k}X_{k}^*.
\end{eqnarray*}
Hence, we get (\ref{FBDE-3}) and (\ref{leader's optimal control-3}).

\emph{ii)$\Rightarrow$i)}. As a fully decoupled FB$\Delta$E, (\ref{FBDE-3}) is solvable. Under the condition, one has
\begin{eqnarray}
\label{X-1}
X_{k+1}^*=(\widetilde{A}_{k}-\widetilde{B}_{k}F_k^{-1}O_k)X_{k}^*+(\widetilde{\mathbf{C}}_{k}^T
-\widetilde{B}_{k}F_k^{-1}\mathbf{D}_k^T)\mathbf{T}_{k+1}X_{k+1}^{*}, ~~~k\in \mathbb{T}_t.
\end{eqnarray}
Under (\ref{T}) and by the first equation of (\ref{FBDE-3}), we have
\begin{eqnarray}
\label{TX-1}
&&\mathbf{T}_{k}X_{k}^{*}=(\mathbf{L}_{k}-\mathbf{K}_{k}F_{k}^{-1}\mathbf{D}_k^T)\mathbf{T}_{k+1}\big{[}\mathbf{I}-(\widetilde{\mathbf{C}}_{k}^T-\widetilde{B}_{k}F_{k}^{-1}\mathbf{D}_k^T)
\mathbf{T}_{k+1}\big{]}^{-1}(\widetilde{A}_{k}-\widetilde{B}_{k}F_{k}^{-1}O_{k})X_k^*  \nonumber \\[1mm]
&&\hphantom{\mathbf{T}_kX_{k}^{*}=}+(\mathbf{H}_{k}-\mathbf{K}_{k}F_{k}^{-1}O_{k})X_k^* \nonumber \\[1mm]
&&\hphantom{\mathbf{T}_kX_{k}^{*}}=(\mathbf{L}_{k}-\mathbf{K}_{k}F_{k}^{-1}\mathbf{D}_k^T)\mathbf{T}_{k+1}X_{k+1}^{*}+(\mathbf{H}_{k}-\mathbf{K}_{k}F_{k}^{-1}O_{k})X_k^*, ~~~k\in \mathbb{T}_t.
\end{eqnarray}
By comparing (\ref{X-1}), (\ref{TX-1}) and (\ref{FBDE-5}), we can see that $(X_k^*, \mathbf{T}_kX_k^*)$ is the solution to (\ref{FBDE-5}).
%
%
By reversing the proof of i)$\Rightarrow$ii), we have that
\begin{eqnarray*}
v_k^*=-F_k^{-1}(O_kX_k^*+\mathbf{D}_k^T\mathbf{Z}_{k+1}^{*}), ~~~k\in \mathbb{T}_t.
\end{eqnarray*}
satisfies the stationary condition (\ref{inequality-2}). Due to the invertibility of $F_k$ and the uniqueness of the solution of (\ref{FBDE-5}), it can be seen that the above $v^*$ is unique. 
The proof is now completed.

\section{Numerical Example}\label{Section-4}
\begin{example}
Consider the Stackelberg game (\ref{cost-1})-(\ref{system-dynamics}) with parameters:
\begin{eqnarray*}
&&A=\left(
\begin{array}{cc}
1& 0.5\\0.3 &2
\end{array}
\right),~~B_1=\left(
\begin{array}{cc}
1& 1\\ 0& 1.2
\end{array}
\right),~~B_2=\left(
\begin{array}{cc}
0.6& 2\\1& 1.6
\end{array}
\right),~~Q_1=\left(
\begin{array}{cc}
1& 0.5\\0.5& 1.5
\end{array}
\right),\\[1mm]
&&Q_2=\left(
\begin{array}{cc}
0.6& 0.2\\0.2& 0.8
\end{array}
\right),~~R_1=\left(
\begin{array}{cc}
0.8& 0.3\\0.3& 1
\end{array}
\right),~~R_2=\left(
\begin{array}{cc}
0& 0\\0& 0
\end{array}
\right),~~W_1=\left(
\begin{array}{cc}
1.25& 0.5\\0.5& 1.4
\end{array}
\right),\\[1mm]
&&W_2=\left(
\begin{array}{cc}
1.45& 0.3\\0.3& 1
\end{array}
\right),~~G_1=\left(
\begin{array}{cc}
1& 0.65\\0.65& 1
\end{array}
\right),~~G_2=\left(
\begin{array}{cc}
0.5& -0.4\\-0.4& 0.5
\end{array}
\right).
\end{eqnarray*}
Letting $t=0$, $N=3$ and $x=(1, 0)^T$, find the open-loop equilibrium solution.
\end{example}

\emph{Solution.} By some calculations, we get the following parameters
\begin{eqnarray*}
&&M_0=\left(
\begin{array}{cc}
2.1841& 2.6175\\
2.6175& 9.1965
\end{array}
\right),~~M_1=\left(
\begin{array}{cc}
2.1360& 2.6922\\
2.6922& 8.5144
\end{array}
\right),~~M_2=\left(
\begin{array}{cc}
1.8000& 2.0800\\
2.0800& 5.0000
\end{array}
\right),~~\\[1mm]
&&F_0=\left(
\begin{array}{cc}
1.4500& 0.3000\\
0.3000& 1.0000
\end{array}
\right),~~F_1=\left(
\begin{array}{cc}
1.4500& 0.3000\\
0.3000& 1.0000
\end{array}
\right),~~F_2=\left(
\begin{array}{cc}
1.4500& 0.3000\\
0.3000& 1.0000
\end{array}
\right),~~\\[1mm]
&&\mathbf{I}-(\widetilde{\mathbf{C}}_{0}^T-\widetilde{B}_{0}F_{0}^{-1}\mathbf{D}_{0}^T)\mathbf{T}_{1}=\left(
\begin{array}{cc}
1.0371& -0.0969\\
0.0417& 0.8908
\end{array}
\right),~~\\[1mm]
&&\mathbf{I}-(\widetilde{\mathbf{C}}_{1}^T-\widetilde{B}_{1}F_{1}^{-1}\mathbf{D}_{1}^T)\mathbf{T}_{2}=\left(
\begin{array}{cc}
0.9552& 0.1173\\
-0.0019& 0.9903
\end{array}
\right),\\[1mm]
&&\mathbf{I}-(\widetilde{\mathbf{C}}_{2}^T-\widetilde{B}_{2}F_{2}^{-1}\mathbf{D}_{2}^T)\mathbf{T}_{3}=\left(
\begin{array}{cc}
0.9971& 0.0650\\
-0.0547& 1.0858
\end{array}
\right).
\end{eqnarray*}
It is easy to see that $M_k$, $k=0,1,2$, are positive define and $F_k$, $\mathbf{I}-(\widetilde{\mathbf{C}}_{k}^T-\widetilde{B}_{k}F_{k}^{-1}\mathbf{D}_{k}^T)\mathbf{T}_{k+1}$, $k=0,1,2$, are invertible. Hence, the problem admits a unique open-loop equilibrium solution for the initial pair $(0, x)$, which is given by
%
%
\begin{eqnarray*}
&&(u^{*(0,x)}_0, v^{*(0,x)}_0)=\left(
\begin{array}{cc}
-0.3711&  0.0053\\
-0.3204& -0.0057
\end{array}
\right),~~~(u^{*(0,x)}_1, v^{*(0,x)}_1)=\left(
\begin{array}{cc}
-0.1583& 0.0230\\
 -0.0632& 0.0462
\end{array}
\right),\\[1mm]
&&(u^{*(0,x)}_2, v^{*(0,x)}_2)=\left(
\begin{array}{cc}
 -0.0456& 0.0254\\
-0.0139& 0.0094
\end{array}
\right).
\end{eqnarray*}
At time instant 1, the equilibrium state $X^{*(0,x)}_1=(0.3003, -0.0883)^T$, which is denoted as $z$. Now, reconsider this Stackelberg game for the initial pair $(1, z)$. Then, the unique open-loop equilibrium solution for the initial pair $(1, z)$ is
\begin{eqnarray*}
(u^{*(1,z)}_1, v^{*(1,z)}_1)=\left(
\begin{array}{cc}
-0.1583& 0.0230\\
 -0.0632& 0.0462
\end{array}
\right),~~~(u^{*(1,z)}_2, v^{*(1,z)}_2)=\left(
\begin{array}{cc}
 -0.0456& 0.0254\\
-0.0139& 0.0094
\end{array}
\right).
\end{eqnarray*}
Clearly,
\begin{eqnarray*}
(u^{*(1,z)}_1, v^{*(1,z)}_1)=(u^{*(0,x)}_1, v^{*(0,x)}_1),~~(u_2^{*(1,z)}, v^{*(1,z)}_2)=(u_2^{*(0,x)}, v^{*(0,x)}_2).
 \end{eqnarray*}
It shows that the open-loop equilibrium solution is time-consistent.  \hfill $\square$

\section{Conclution}\label{Section-5}

In this paper, open-loop equilibrium solution is investigated for deterministic dynamic Stackelberg game, which is shown to be time-consistent. 
Necessary and sufficient condition for the existence and uniqueness of open-loop equilibrium solution is given, and two Riccati-like equations are introduced to characterize the open-loop equilibrium solution. For future research, we may study the open-loop equilibrium solution for stochastic Stackelberg games.

%

\end{document}